\documentclass[12pt,a4paper]{amsart}
\usepackage{amsmath, amssymb, amsthm,
amsfonts, amsthm, bm, eucal, graphicx}
\makeatletter
\@namedef{subjclassname@2020}{\textup{2020} Mathematics Subject Classification}
\makeatother
\usepackage[all]{xy}
\usepackage{graphics}
\usepackage{epsfig}
\usepackage{graphicx}
\theoremstyle{plain}
\usepackage{amssymb}
\usepackage[colorlinks]{hyperref}
\usepackage{enumerate}

\advance\hoffset-20mm \advance\textwidth40mm

\newtheorem{thm}{Theorem}[section]
\newtheorem{cor}[thm]{Corollary}
\newtheorem{lem}[thm]{Lemma}
\newtheorem{prop}[thm]{Proposition}
\newtheorem{que}[thm]{Question}

\theoremstyle{definition}
\newtheorem{defi}[thm]{Definition}
\newtheorem{defis}[thm]{Definitions}
\newtheorem{conj}[thm]{Conjecture}
\newtheorem{conv}[thm]{Convention}
\newtheorem{nota}[thm]{Notation}
\newtheorem{rem}[thm]{Remark}
\newtheorem{rems}[thm]{Remarks}
\newtheorem{exa}[thm]{Example}
\newtheorem{exas}[thm]{Examples}
\newtheorem{sit}[thm]{}

\newcommand{\brem}{\begin{rem}}
\newcommand{\brems}{\begin{rems}}
\newcommand{\erem}{\end{rem}}
\newcommand{\erems}{\end{rems}}
\newcommand{\bexa}{\begin{exa}}
\newcommand{\bexas}{\begin{exas}}
\newcommand{\eexa}{\end{exa}}
\newcommand{\eexas}{\end{exas}}
\newcommand{\bdefi}{\begin{defi}}
\newcommand{\edefi}{\end{defi}}
\newcommand{\bdefis}{\begin{defis}}
\newcommand{\edefis}{\end{defis}}
\newcommand{\bcor}{\begin{cor}}
\newcommand{\ecor}{\end{cor}}
\newcommand{\blem}{\begin{lem}}
\newcommand{\elem}{\end{lem}}
\newcommand{\bconv}{\begin{conv}}
\newcommand{\econv}{\end{conv}}
\newcommand{\bconj}{\begin{conj}}
\newcommand{\econj}{\end{conj}}
\newcommand{\bprop}{\begin{prop}}
\newcommand{\eprop}{\end{prop}}
\newcommand{\bthm}{\begin{thm}}
\newcommand{\ethm}{\end{thm}}
\newcommand{\bnota}{\begin{nota}}
\newcommand{\enota}{\end{nota}}
\newcommand{\bsit}{\begin{sit}}
\newcommand{\esit}{\end{sit}}
\newcommand{\be}{\begin{equation}}
\newcommand{\ee}{\end{equation}}
\newcommand{\bproof}{\begin{proof}}
\newcommand{\eproof}{\end{proof}}
\newcommand{\bque}{\begin{que}}
\newcommand{\eque}{\end{que}}
\def\ba{\begin{array}}
\def\ea{\end{array}}
\def\bea{\begin{eqnarray}}
\def\eea{\end{eqnarray}}

\def\bnum{\begin{enumerate}}
\def\enum{\end{enumerate}}
\newcommand{\la}{\label}

\newtheorem*{theo*}{Theorem}

\theoremstyle{definition}

\newtheorem*{definition*}{Definition}

\def\fg{{\mathfrak g}}
\def\fm{{\mathfrak m}}
\def\fl{{\mathfrak l}}

\def\cH{{\mathcal H}}

\def\cN{{\mathcal N}}
\def\cO{{\mathcal O}}

\def\AA{{\mathbb A}}

\def\NN{{\mathbb N}}
\def\ZZ{{\mathbb Z}}
\def\QQ{{\mathbb Q}}
\def\TT{{\mathbb T}}
\def\kk{{\Bbbk}}
\def\CC{{\mathbb C}}

\def\Ga{{\mathbb G_{\mathrm a}}}

\def\RR{{\mathbb R}}

\def\fb{{\mathfrak b}}
\def\fh{{\mathfrak h}}

\def\fj{{\mathfrak j}}
\def\fl{{\mathfrak l}}

\def\ft{{\mathfrak t}}
\def\fu{{\mathfrak u}}
\def\fb{{\mathfrak b}}

\def\span{\mathop{\rm span}}
\def\End{\mathop{\rm End}}
\def\Aut{\mathop{\rm Aut}}
\def\Lie{\mathop{\rm Lie}}

\def\SAut{\mathop{\rm SAut}}

\def\JONQ{\mathop{\rm JONQ}}

\def\Ad{\mathop{\rm Ad}}
\def\ad{\mathop{\rm ad}}

\def\supp{\mathop{\rm supp}}

\def\Cent{\mathop{\rm Cent}}

\def\Spec{\mathop{\rm Spec}}

\def\Der{\mathop{\rm Der}}
\def\Vec{\mathop{\rm Vec}}

\def\deg{\mathop{\rm deg}}

\def\p{\partial}

\def\rk{\mathop{\rm rk}}

\def\Aff{\mathop{\rm Aff}}

\def\ll1{l_{\lambda}^{-1}(1)}
\def\lm1{l_{\mu}^{-1}(1)}

\def\ba{\begin{array}}
\def\ea{\end{array}}
\def\bea{\begin{eqnarray}}
\def\eea{\end{eqnarray}}

\begin{document}
\sloppy
\title[Lie algebras of derivations]
{Locally finite solvable Lie algebras of derivations}

\author{Mikhail\ Zaidenberg}
\address{Univ. Grenoble Alpes, CNRS, IF, 
38000 Grenoble, France}
\email{mikhail.zaidenberg@univ-grenoble-alpes.fr}
%
%
\begin{abstract} 
Let $X$ be an affine variety. 
The local finiteness of 
a Lie subalgebra
$\fh\subset\Lie(\Aut(X))$ is equivalent to 
the existence of an algebraic subgroup 
$G\subset\Aut(X)$ such that 
$\fh\subset\Lie(G)$.
Let $\fh$ be a solvable Lie subalgebra  of 
$\Lie(\Aut(X))$
generated by a finite collection of
locally finite Lie subalgebras. The authors 
of \cite{AZ25b} wondered
whether 
$\fh$ is itself locally finite. 
After presenting
some criteria for the local finiteness of $\fh$,
we answer this question in the  affirmative  
in the particular case where $X$ 
is the affine plane $\AA^2$.
\end{abstract}
\subjclass[2020]{Primary 14J50, 17B66; \ 
Secondary 14R10, 17B45, 22E65}
\keywords{Affine varieties, 
automorphism group, Lie algebra, 
solvable subalgebra,
derivation, 
vector field, locally finite}
\maketitle
{\footnotesize \tableofcontents}

\section{Introduction}
All algebraic 
varieties in this paper are defined 
over an algebraically closed field $\kk$
of characteristic zero. 
Throughout the article,
$X$ denotes an affine variety over $\kk$,
and $\fh$ a Lie subalgebra of
$\Lie(\Aut(X))$. We say that $\fh$
is \emph{locally finite} if every $f\in\cO(X)$
belongs to a finite-dimensional 
vector subspace 
of $\cO(X)$ invariant under $\fh$,
see, e.g., \cite[Definition~1.5.1]{KZ24}.
It is known that 
a subalgebra $\fh$ of $\Lie(\Aut(X))$ 
is locally finite if and only if $\fh\subset\Lie(G)$,
where $G\subset\Aut(X)$ is 
a connected algebraic 
subgroup, see 
\cite[Theorem~1]{CD03} 
or \cite[Theorem~E]{KZ24}.

A derivation $\p\in\Der(\cO(X))$ is
said to be locally finite if the Lie subalgebra
$\kk\p$ is locally finite. Thus, any  locally finite
Lie subalgebra $\fh\subset\Lie(\Aut(X))$
consists of locally finite derivations, and is, 
consequently, generated by 
locally finite derivations. 

A Lie subalgebra $\fh\subset\Lie(\Aut(X))$
is said to be \emph{locally integrable}
(following the terminology of  \cite{AZ25b}),
or \emph{weakly locally finite}
(following the terminology of \cite{CDR26})
if $\fh$ is filtered by an 
increasing sequence $\fh_i$ of 
locally finite Lie subalgebras. 
If  $\fh$ is weakly locally finite,
then every finitely generated 
Lie subalgebra of $\fh$ is contained in 
some $\fh_i$, and is therefore 
locally finite. 

Recall the following question.
\bque[{\rm \cite[Question~1]
{AZ25b}}]\la{ques-1}
Let $\fh=\langle \fh_1,
\ldots,\fh_k\rangle_{\Lie}$
(resp. 
$\fh=\langle \fh_1,\ldots,\fh_k
\ldots\rangle_{\Lie}$)
be a solvable Lie subalgebra of 
$\Lie(\Aut(X))$ generated by 
locally finite Lie subalgebras 
$\fh_i\subset \Lie(\Aut(X))$, 
$i=1,2,\ldots$. 
Is it true that $\fh$ is locally finite
(resp. weakly locally finite)?
\eque
The author knows of no example where 
a finitely generated $\fh$ 
as in Question \ref{ques-1} 
is not locally finite. 
We therefore conjecture that 
Question  \ref{ques-1} 
has an affirmative answer. 
The following theorem 
confirms this conjecture 
in the case where $X$ 
is the affine plane $\AA^2=\AA^2_{\kk}$.
\bthm\la{main-thm}
A solvable Lie subalgebra 
\[ \fh=\langle \p_1,\ldots,\p_k
\ldots\rangle_{\Lie} 
\subset\Lie(\Aut(\AA^2))\]
generated by locally
finite derivations 
$\p_i$ is triangulable  and 
weakly locally finite.
If $\fh$ is generated by 
a finite set of 
locally finite derivations, 
then it is locally finite.
\ethm
Recall that
a Lie subalgebra 
$\fh\subset\Der(\kk[x,y])$ 
is said to be \emph{triangulable}
if $\fh$ is $\Ad$-conjugate to a 
Lie subalgebra consisting 
of triangular derivations, 
see subsection \ref{ss:3-1}. 

Note that a locally finite Lie subalgebra 
$\fh\subset\Lie(\Aut(X))$
is finite-dimensional, 
see \cite[Lemma~1.6.2]{KZ24}. 
Question \ref{ques-1}
 is related to the following one.
\bque[{\rm  see \cite[Question~2]{KZ24}}]
Let $\{\fh_i\}_i$ be a family 
of locally finite Lie 
subalgebras of $\Lie(\Aut(X))$.
Is it true that the Lie subalgebra 
$\fh\subset\Lie(\Aut(X))$ 
generated by the $\fh_i$
is locally finite provided 
that it is finite-dimensional?
\eque
See  \cite[Sec. 1.6]{KZ24}  
and Corollary \ref{cor:lf-fd} 
below for some partial results.
Here, we only consider the case of 
a solvable Lie algebra $\fh$. 

Recall that a derivation $\p\in \Der(\cO(X))$ 
is said to be \emph{locally nilpotent} 
if every regular function $f\in\cO(X)$
is annihilated by an iterate $\p^{(n)}$, 
where $n=n(f)$. Every locally nilpotent
derivation is locally finite. 
If $\fh$ as in Theorem \ref{main-thm}
is generated by
locally nilpotent derivations, 
then it 
consists of locally nilpotent 
derivations. 
Any finitely generated subalgebra 
$\fh'$ of $\fh$
coincides with  the Lie 
algebra of a unipotent algebraic 
subgroup of $\Aut(X)$. 
Consequently, $\fh'$ is 
locally finite, and $\fh$ 
is weakly locally finite,
see Theorem \ref{thm:unip}.

Let us review some known related 
results. 
A Lie subalgebra of 
$\Lie(\Aut(\AA^2))$ 
\underline{consisting of} locally 
nilpotent derivations is
triangulable, see
\cite[Theorem~3.11]{PS17}.  
Consequently, it is 
solvable, locally nilpotent, and
weakly locally finite. 
The local nilpotence and 
weak local finiteness survive
for Lie subalgebras of $\Der(X)$
\underline{consisting of} locally 
nilpotent derivations, where $X$
is a general, not necessarily 
reduced affine variety, see
\cite[Theorem~1.4]{BPZ24}
and \cite{Sku21}.
 
Our Theorem \ref{main-thm} 
 implies the triangulability of
 a solvable Lie subalgebra $\fh$ of 
$\Lie(\Aut(\AA^2))$ consisting
of locally finite derivations, 
thereby answering, 
in this particular setting, 
a question posed 
to the author by Andriy Regeta. 

More generally, in
the recent article \cite{CDR26},
it is established that,
for any quasi-affine variety $Y$ 
defined over an arbitrary field,
any solvable
Lie subalgebra of $\Lie(\Aut(Y))$
\underline{consisting of} locally finite
derivations 
is weakly locally finite 
(see \cite[Corollary~4]{CDR26}).

Another special case
where triangulability holds
is the following; see \cite{Sku21}
(cf. also \cite[Theorem~D]{KZ24}). 
Let $\fh$
be a Lie subalgebra of 
$\Lie(\Aut(\AA^n))$ 
\underline{consisting
of} locally nilpotent 
derivations. Suppose that 
the intersection of their 
kernels consists only of constants.
Then $\fh$ is triangulable, 
and consequently solvable 
and weakly locally finite. 
If, moreover, $\fh$ is maximal
among Lie subalgebras consisting 
of locally nilpotent derivations, 
then $\fh$ coincides,
in a suitable
coordinate system on $\AA^n$, 
with the Lie algebra of triangular 
locally nilpotent  derivations; 
see \cite[Theorem~1]{Sku21}. 
If the aforementioned condition 
regarding
the kernels is not satisfied, then 
the triangulability may also fail to hold,
as is the case in the example 
due to Bass \cite{Bas84}; see 
Remark \ref{rem:Bass}. 
The solvability
of a Lie subalgebra 
$\fh\subset\Lie(\Aut(X))$
\underline{consisting of} 
locally nilpotent derivations 
is established
in \cite[Corollary~3.4]{Per23} 
for any affine variety $X$.

We believe that an analogue 
of Theorem \ref{main-thm}, 
with a similar proof, 
also remains valid for 
normal affine toric surfaces.  
As Bass' example shows, 
the triangulability assertion
in Theorem \ref{main-thm} 
is not valid, in general, for 
solvable locally finite Lie subalgebras 
of $\Lie(\Aut(\AA^n))$
in higher dimensions, see 
Remark \ref{rem:Bass}. 

The content of the paper is as follows.
In Section \ref{sec-2} we consider
Question \ref{ques-1} for general 
affine varieties, and we perform 
several successive reductions.

In Section \ref{sec-3},
we address the case of 
the affine plane. 
After recalling in subsections 
\ref{ss:3-1}--\ref{ss:comm-lnd}
the necessary 
preliminaries (on 
the structure of 
$\Lie(\Aut(\AA^2))$ as 
a bigraded Lie algebra, 
on the technique of 
Newton polygons, 
etc.), 
we prove, 
in subsections 
\ref{ss:solv}--\ref{ss:3.8},
Theorem \ref{main-thm} 
under different 
additional assumptions.
First, we assume that 
the Lie subalgebra $\fh$
from Question \ref{ques-1} 
is generated 
by locally nilpotent derivations 
(see Theorem \ref{thm:unip}), 
and then we treat the general case
in Theorem \ref{thm:rk-1}. 
The proof of 
Theorem \ref{main-thm}
 is completed at the end 
 of Section \ref{sec-3}.
\section{General affine 
varieties}\la{sec-2}
\subsection{First reduction}
We can assume that every $\fh_i$ 
in Question \ref{ques-1} is 
one-dimensional. 
Indeed, being locally finite, $\fh_i$ 
is finite-dimensional
and consists of locally finite derivations. 
Thus, every $\fh_i$ is generated by 
a finite set of locally finite derivations. 
\subsection{Second reduction} 
Recall that a derivation 
$\delta\in\Der(\cO(X))$ 
is said to be \emph{semisimple} 
if $\delta$ is locally finite
and if the restriction of $\delta$ to 
any finite-dimensional $\delta$-invariant 
vector subspace $V\subset \cO(X)$ 
is diagonalizable. 
Let us fix an ascending 
filtration of $\cO(X)$ by 
a sequence $V_n$ of such subspaces, 
as well as the finite decompositions 
$V_n=\oplus_{\lambda\in\kk} 
V_{n, \lambda}$,
where $V_{n, \lambda}\subset V_n$ is the 
eigenspace of $\delta|_{V_n}$
that corresponds to the eigenvalue 
$\lambda$. It is clear that
$V_{n, \lambda}\subset 
V_{n+1, \lambda}$ for all $n\ge 1$ 
and all $\lambda\in\kk$. 
A basis of $V_{n,\lambda}$
consisting of eigenvectors of $\delta$
can be extended to such a basis
of $V_{n+1, \lambda}$;
this inductive process yields
a countable basis of $\cO(X)$
consisting of eigenvectors of $\delta$.
The restriction of $\delta$
to any $\delta$-invariant subspace 
$V\subset\cO(X)$
is also semisimple. 

Every locally finite 
derivation $\p\in\Lie(\Aut(X))$ 
admits a unique Jordan decomposition
$\p=\p_{\mathrm s}+\p_{\mathrm n}$, 
where $\p_{\mathrm s}$ is semisimple 
and $\p_{\mathrm n}$ is locally nilpotent;
they commute, 
and every vector subspace 
$V\subset\cO(X)$
invariant under $\p$
is  also invariant under 
$\p_{\mathrm s}$ and $\p_{\mathrm n}$, 
see, e.g.,  
\cite[Theorems~9.4.2 and 9.4.3]{Nov94}.
Furthermore, 
$\ad_{\p}=\ad_{\p_{\mathrm s}}
+\ad_{\p_{\mathrm n}}$ is 
a Jordan decomposition of $\ad_{\p}$
acting on $\Lie(\Aut(X))$ 
(see \cite[Lemma~1.7.1]{KZ24}), and 
any vector subspace 
$V\subset \Lie(\Aut(X))$ 
invariant under $\ad_{\p}$ 
 is also invariant 
under $\ad_{\p_{\mathrm s}}$ 
and $\ad_{\p_{\mathrm n}}$. 
This follows from 
\cite[Theorem~9.4.2(3)]{Nov94}
combined with 
\cite[Lemma~1.7.1]{KZ24}.

These observations lead to 
the following lemma.
For a Lie algebra $\fh$, we set, 
as usual, $\fh^{(1)}=[\fh,\fh]$ 
and $\fh^{(k+1)}=[\fh^{(k)},\fh^{(k)}]$. 
\blem\la{lem:red-2}
For a solvable Lie subalgebra 
 $\fh\subset\Lie(\Aut(X))$ 
of derived length $l$,
the following is true.
\begin{itemize}
\item[$(a)$] Let
$\p\in\fh$ be 
a locally finite derivation.
Then the Lie algebra 
$\widehat{\fh}=\langle \fh, 
\p_{\mathrm s}\rangle_{\Lie}$
is solvable 
of derived length $\le l+1$. 
\item[$(b)$] 
Let  $\fg\subset\Lie(\Aut(X))$ 
be an abelian Lie subalgebra
that normalizes $\fh$, that is, 
$\ad_{\fg}(\fh)\subset\fh$. 
Then the Lie algebra 
$\widehat{\fh}=\langle \fh, 
\fg\rangle_{\Lie}$ is solvable 
of derived length $\le l+1$. 
\end{itemize}
\elem
\bproof
(a) Note that $\fh$ 
is invariant under 
$\ad_{\p}=\ad_{\p_{\mathrm s}}
+\ad_{\p_{\mathrm n}}$
and 
any $\ad_{\p}$-invariant subspace is stable 
under both summands. 
In particular, $\fh$ is invariant under 
$\ad_{\p_{\mathrm s}}$. 
Then
\[\fh^{(1)}=[\widehat{\fh},\widehat{\fh}]
=[\fh,\fh]+[\p_s,\fh]\subset \fh,\]
since $[\p_{\mathrm s},\p_{\mathrm s}]=0$.
From this one obtains inductively
\be\la{eq:h}\widehat{\fh}^{(n+1)}\subset
\fh^{(n)}.\ee
Now statement (a) follows. 

 (b) Under the assumptions of (b), we have 
$\ad_{\fg}(\fh)\subset\fh$ and $[\fg,\fg]=0$. 
It follows, as before, 
that $\widehat{\fh}^{(1)}\subset\fh$.
Consequently, \eqref{eq:h} holds again, 
which proves (b). 
\eproof
Recall that a Lie subalgebra
$\fh\subset\Lie(\Aut(X))$ 
is said to be
\emph{$J$-saturated}, if, for every locally 
finite derivation
$\p\in\fh$, it
contains its locally nilpotent 
and semisimple parts
$\p_{\mathrm n}$ and $\p_{\mathrm s}$
(see \cite[Section~3.2]{KZ24}). 
By virtue of the following proposition, 
when addressing Question \ref{ques-1}
one may assume that $\fh$ 
is maximal (in a certain sense) 
and $J$-saturated. 
\bprop\la{prop:max} 
Let $\fh\subset\Lie(\Aut(X))$ 
be a solvable Lie subalgebra 
generated by locally finite derivations.
Then $\fh$ is contained in a 
Lie subalgebra 
$\fh_{\max}\subset\Lie(\Aut(X))$ 
that is solvable, generated by 
locally finite derivations, 
and maximal among the
Lie subalgebras
of $\Lie(\Aut(X))$ possessing these 
two properties. Any Lie subalgebra
$\fh_{\max}$ of this type
is \emph{$J$-saturated}.
\eprop
\bproof Let $\cH$ be the family of all 
solvable Lie subalgebras  of $\Lie(\Aut(X))$
generated by locally finite derivations.
Recall that the length 
of the derived series of
any solvable Lie subalgebra 
of $\Der(\cO(X))$
is at most $2\dim(X)$
(see \cite[Corollary~4]{MP14}). 
It follows that the union of 
any increasing chain 
$(\fh_\alpha)$ of 
elements of $\cH$ belongs to $\cH$.
By Zorn's lemma, 
the first assertion follows.

Let $\p=
\p_{\mathrm n}+\p_{\mathrm s}$
be the Jordan decomposition
 of a locally finite derivation 
$\p\in\fh_{\max}$. 
By Lemma \ref{lem:red-2}, 
the Lie algebra 
$\widehat{\fh_{\max}}=\langle \fh_{\max},
\p_{\mathrm s}\rangle_{\Lie}$
also belongs to $\cH$. Consequently,
$\widehat{\fh_{\max}}=\fh_{\max}$.
Hence, $\fh_{\max}$ contains 
$\p_{\mathrm s}$ and $\p_{\mathrm n}$;
it is therefore $J$-saturated. 
\eproof
The following corollary is immediate.
\bcor\la{cor:red-2}
The answer to 
Question \ref{ques-1} is affirmative
if and only if it is affirmative 
under the additional 
assumption that 
$\fh$ in Question \ref{ques-1} is 
generated by 
locally nilpotent 
and semisimple 
derivations. 
\ecor
\subsection{Third reduction} 
We use the following terminology, 
see for example \cite{Hoch81} and 
\cite[Sec.~1.8]{KZ24}.
\bdefis  $\,$

1. 
A Lie subalgebra 
$\fh\subset\Der(\cO(X))$
 is said to be \emph{algebraic} 
(resp. \emph{integrable})
if $\fh=\Lie(G)$ (resp. 
$\fh\subset \Lie(G)$) 
for an algebraic subgroup 
$G\subset\Aut(X)$. 
A derivation $\p\in\Der(\cO(X))$ 
is said to be 
\emph{algebraic} 
(resp. \emph{integrable})
if the Lie subalgebra 
$\kk\p$ is algebraic 
(resp. integrable). 

2. A Lie subalgebra 
$\ft\subset\Lie(\Aut(X))$ 
is said to be \emph{toral} if it consists 
of semisimple derivations. The \emph{rank} 
$\rk(\fh)$ is the maximal dimension 
of toral subalgebras of $\fh$.
\edefis
We know that $\fh$ is integrable 
if and only if it is locally finite, 
see \cite[Theorem~1]{CD03}
or \cite[Theorem~E]{KZ24}.
In the latter case, $\fh\subset\Lie(\Aut(X))$. 
For example, every locally nilpotent derivation 
is algebraic, and every semisimple 
derivation is integrable. Indeed, 
if $s\in\Der(\cO(X))$ is semisimple,
then $s$ is contained in a unique minimal 
algebraic toral subalgebra
$\ft_{\min}(s)\subset\Lie(\Aut(X))$, 
see the following lemma.
 \blem[{\rm \cite[Lemma~1.8.2]{KZ24}}] 
 \la{lem:toral} 
 For a toral Lie subalgebra 
 $\ft \subset\Der(\cO(X))$, 
 there exists a unique
 smallest torus 
 $T_{\min}(\ft) 
 \subset \Aut(X)$ such that 
 $\ft\subset\ft_{\min}(\ft):= 
 \Lie(T_{\min}(\ft))$. 
Every vector subspace 
$E \subset\Der(\cO(X))$ 
invariant under $\ad_{\ft}$ 
 is also invariant under $T_{\min}(\ft)$.
Furthermore, $\ft$ is locally finite and 
 $\dim(\ft)\le\dim(X)$.
If $\dim(\ft)=\dim(X)$, 
 then $X$ is a toric variety and 
 $\ft$ is algebraic.
 \elem
In the next Propositions 
 \ref{prop:integr}--\ref{prop:KZ} 
 and Corollary \ref{cor:lf-fd} we
 provide criteria for  local finiteness.
 \bprop\la{prop:integr}
Let $\fh \subset \Lie(\Aut(X))$ be 
a solvable Lie subalgebra. 
Then $\fh$ is locally finite
if and only if there exists 
a solvable connected 
algebraic subgroup 
$G \subset \Aut(X)$ 
with derived length 
$\le \dim(X) + 1$ such that 
$\fh \subset \Lie(G)$.
Such a minimal
algebraic subgroup 
$G_{\min}(\fh)$ is unique. 
\eprop
\bproof 
See  \cite[Example~1.6.1]{KZ24}
for the ``if'' direction in the first 
statement and 
\cite[Lemma~2.21]{AZ25b}
for the``only if''direction. 
The uniqueness statement 
follows from the minimality condition. 
\eproof
\bprop[{\rm \cite[Sec. 1.6, 
Corollary]{KZ24}}] 
\la{prop:KZ}
Let $\fh\subset \Lie(\Aut(X))$ 
be the Lie subalgebra 
generated by a family of locally 
finite Lie subalgebras 
$\fh_i \subset\Lie(\Aut(X))$, $i\in I$. 
 Suppose that $\fh_i$ 
 is algebraic  for all $i\in I$, that is,
 $\fh_i  = \Lie(G_i)$ for a connected 
 algebraic group 
 $G_i\subset \Aut(X)$. 
 Then $\fh$ is locally finite 
 if and only if $\fh$ 
 is finite-dimensional. In the latter case, 
 the subgroup 
 $G$ generated by the $G_i$ 
 is algebraic and 
 $\fh= \Lie(G)$. 
\eprop
\bcor\la{cor:lf-fd} 
Consider a Lie subalgebra
\be\la{eq:gener}
\fh=\langle a_1,\ldots, 
a_k, b_1,\ldots,b_l\rangle_{\Lie}
\subset \Lie(\Aut(X)),\ee
where 
$a_1,\ldots, a_k$ are locally nilpotent and
$b_1,\ldots,b_l$ are semisimple.
Suppose that all generators $b_i$ 
are algebraic,
that is, for all $i=1,\ldots, l$ we have 
$\kk b_i=\Lie(T_i)$ for 
a one-dimensional algebraic 
torus $T_i\subset\Aut(X)$. Then $\fh$ 
is locally finite if and only if it 
is finite-dimensional.
\ecor 
\bproof For $i=1,\ldots,k$ 
we have $\kk a_i=\Lie(U_i)$,
where 
$U_i=\{\exp(ta_i)\,|\,t\in\kk\}$
is a one-parameter 
unipotent subgroup of $\Aut(X)$. 
The assertion then follows
from Proposition \ref{prop:KZ}.
\eproof
\bque\la{ques-2} Let 
$\fh\subset \Lie(\Aut(X))$
be a Lie subalgebra  
generated by a finite family of 
locally nilpotent and semisimple 
derivations as in  \eqref{eq:gener}, 
where, this time, we do not 
assume that all 
the semisimple generators are algebraic. 
Suppose that $\fh$ is finite-dimensional 
(resp. solvable),
and consider the Lie subalgebra 
\be\la{eq:tilde-h}\widetilde\fh=\langle 
\kk a_1,\ldots,\kk a_k, \Lie(T_{\min}(b_1)),
\ldots, \Lie(T_{\min}(b_l))\rangle_{\Lie}.\ee
Is $\widetilde\fh$ also finite-dimensional 
(resp. solvable)?
\eque
An affirmative answer 
for the solvable case follows from 
Corollary \ref{cor:alg-gener} below. 
According to
Proposition \ref{prop:KZ}, a positive answer 
for the finite-dimensional case 
would allow us to answer 
Question \ref{ques-1} in the affirmative. 
\brem\la{rem:G} We use the notation from 
Question \ref{ques-2}.
Consider the subgroup 
\be\la{eq:G} G=\langle U_i=\exp(\kk a_i), T_{\min}(b_j)\,|\,
i=1,\ldots,k,\,j=1,\ldots, l\rangle 
\subset\Aut(X).\ee
Suppose that $G$ is solvable. Then 
$G$ is a connected algebraic subgroup,
see \cite[Theorem~A]{CKRvS26}.
 Consequently, 
the Lie algebra
$\widetilde \fh\subset\Lie(G)$ 
is solvable, finite-dimensional, 
and algebraic; 
it is therefore locally finite. 

However, in general, 
we do not know whether 
$G$ is solvable. If $G$ is not solvable, then 
$\Lie(G)$ does not need to be solvable. 
However, we do not know
whether $\tilde \fh=\Lie(G)$; apriori, the inclusion 
$\tilde h\subset\Lie(G)$ can be strict, 
and in this case, $\tilde \fh$
can be solvable even if $\Lie(G)$ is not. 
\erem
\subsection{Maximal solvable subalgebras 
generated by locally finite derivations}
The following observations 
will be useful later on.
\blem\la{lem:ideal}
Let $\fh\subset\Lie(\Aut(X))$
be a $J$-saturated 
solvable Lie subalgebra 
generated by 
locally finite derivations. 
Consider the Lie subalgebra 
$\fh_{\rm lnd}$
of $\fh$ generated by 
all locally nilpotent 
derivations of $\fh$. 
Then $\fh_{\rm lnd}$
is an ideal of $\fh$. 
\elem
\bproof 
Being $J$-saturated 
and of countable dimension, 
$\fh$ is generated by 
two sequences, one consisting 
of locally nilpotent
derivations and the other of 
semisimple derivations. 
The elements of the first sequence 
belong to $\fh_{\rm lnd}$. 
It therefore suffices to verify
that $\fh_{\rm lnd}$ is stable
under $\ad_\delta$ for any
semisimple derivation $\delta\in\fh$.
This holds if $\fh_{\rm lnd}$ 
is stable under $\Ad_T$, 
where $T=T_{\min}(\delta)$, 
see \cite[Lemma~1.8.2 
and its proof]{KZ24}. 
Note that $\Ad_T$ 
leaves $\fh$ invariant 
because $\delta$ does, 
see Lemma \ref{lem:toral}. 
Recall that any conjugate of 
 a unipotent subgroup
 of $\Aut(X)$ is unipotent. 
Consequently, $\Ad_T$
transforms locally nilpotent derivations
into locally nilpotent derivations.
Thus, $\Ad_T$ preserves
the Lie subalgebra $\fh_{\rm lnd}$.
It follows that 
$\ad_\delta\in\ad_{\Lie(T)}$ also 
preserves $\fh_{\rm lnd}$.
The proof is thus  complete.
\eproof
\bprop\la{prop:alg-toral}
Let 
$\fh_{\max}\subset\Lie(\Aut(X))$
be a solvable Lie algebra generated 
by locally finite derivations
that is
maximal within the class of 
all Lie subalgebras of this type
(see Proposition \ref{prop:max}). 
Then any maximal toral 
subalgebra $\ft_{\max}$ 
of $\fh_{\max}$ is algebraic. 
\eprop
\bproof
By Lemma \ref{lem:toral}, 
$\ft_{\max}$ is contained in a unique
smallest algebraic toral Lie subalgebra
$\ft_{\rm alg}=\ft_{\rm alg}(\ft_{\max}):
=\Lie(T_{\min}(\ft_{\max}))$, which
normalizes $\fh_{\max}$. 
According to Lemma \ref{lem:red-2}(b),
the Lie subalgebra $\widehat{\fh} =
\langle\fh_{\max},\ft_{\rm alg}\rangle_{\Lie}$
of $\Lie(\Aut(X))$ is solvable and, 
by construction, 
is generated by 
locally finite derivations. 
Therefore, $\widehat{\fh}= 
\fh_{\max}$. 
In particular, 
$\ft_{\max}\subset \ft_{\rm alg}$
coincides with 
$\ft_{\rm alg}\subset\fh_{\max}$.
The  assertion follows. 
\eproof
\bcor\la{cor:alg-gener}
Any solvable Lie subalgebra 
$\fh\subset\Lie(\Aut(X))$
generated by locally finite derivations
is contained in a solvable Lie subalgebra 
$\fh_{\max}\subset\Lie(\Aut(X))$
generated by algebraic Lie subalgebras.
\ecor
\bproof
The assertion follows
immediately from 
Proposition \ref{prop:alg-toral}.
\eproof
\brems\la{rem:filtr} 1. Consequently, 
in Question \ref{ques-1} 
we may suppose that every $\fh_i$
is algebraic. 

2. In the course of the proof of 
Theorem \ref{main-thm},
we obtain an additional result 
valid for arbitrary
affine varieties; see 
Lemma \ref{lem:1-lnd}.
\erems
\section{Case of the 
affine plane}\la{sec-3}
 \subsection{Triangular subalgebras of 
$\Lie(\Aut(\AA^2))$}\la{ss:3-1}
We use the following notation.
\bnota\la{not:u-plus}
 We define
\be\la{eq:u2-plus} \fu^+_2=\kk[y]\frac{\p}{\p x}\oplus\kk\frac{\p}{\p y}
\quad\left(\text{resp.}\quad
\fu^-_2=\kk\frac{\p}{\p x}\oplus\kk[x]\frac{\p}{\p y}\right).\ee
It is well known that
$\fu_2^+$ (resp. $\fu^-_2$) is a 
metabelian (i.e., 2-step solvable) 
subalgebra of 
$\Lie(\Aut(\AA^2))$ of infinite dimension. 
It consists of locally nilpotent 
upper (resp. lower)
triangular derivations. 
The Lie algebra
$\fu_2^+$ is filtered by 
an increasing sequence of locally finite
Lie subalgebras 
\be\la{eq:fu-2} 
\fu_{2,\,\le d}^+=\biggl\{\p=p(y)\frac{\p}{\p x}+c\frac{\p}{\p y}\in\fu^+_2\,|\,
p\in\kk[y],\,\deg(p)\le d,\,\,c\in\kk\biggr\}=\Lie(U_{2,\,\le d}^+),\ee
where $U_{2,\,\le d}^+=
\exp(\fu_{2,\,\le d}^+)$. 
The group $U_2^+=\exp(\fu_2^+)$ is 
an infinite-dimensional  nested unipotent
subgroup of $\Aut(\AA^2)$ filtered by the
unipotent algebraic subgroups 
$(U_2^+)_{\le d}$. 

Consider also the maximal toral 
Lie subalgebra 
\[\ft_2=\kk x\frac{\p}{\p x}\oplus
\kk y\frac{\p}{\p y}\subset
\Der(\kk[x,y]).\]
The semidirect product 
$\fj_2^\pm=\fu^\pm_2\rtimes \ft_2$
is the Lie algebra of the 
de Jonqui\`eres group
$\JONQ^\pm(\AA^2)$ of upper 
(resp. lower) triangular
automorphisms of $\AA^2$.
\enota
\brems\la{rem:fj-2} 
1. Note that $\fj_2^+$
 is solvable
with derived length three; 
indeed, by \eqref{eq:u2-plus},
\[[\fj_2^+,\fj_2^+]=\fu_2^+,\qquad
[\fu_2^+,\fu_2^+]=
\kk[y]\frac{\partial}{\partial x},
\quad\text{and}\quad
\left[\kk[y]\frac{\p}{\p x},
\kk[y]\frac{\p}{\p x}\right]=0.\]
Moreover, $\fj_2^+$ is filtered by the 
increasing sequence of locally finite 
Lie subalgebras 
\[\fj_{2,\,\le d}^+:=
\fu_{2,\,\le d}^+\rtimes \ft_2.\]
Consequently, it is weakly locally finite.

2. Being solvable
with derived length two, $\fu_2^+$ 
is not nilpotent,
and the nilpotency class of 
a locally finite Lie subalgebra
$\fh\subset \fj_2^+$ can be arbitrarily large.
Indeed,
$[p(y)\frac{\p}{\p x},\frac{\p}{\p y}]
=-p'(y)\frac{\p}{\p x},$
so that iterated commutators 
lower the degree 
but never vanish uniformly.
Similarly, the finite-dimensional truncations 
$\fu_{2,\,\le d}^+$
can have arbitrarily large nilpotency class.
\erems
\bdefis\la{def:triangulable} $\,$

1. A Lie subalgebra 
$\fh\subset\Der(\kk[x,y])$ is said to be 
\emph{upper} (resp.  \emph{lower}) 
\emph{triangular}
if $\fh\subset\fj_2^+$ 
(resp. $\fh\subset\fj_2^-$). 
It is said to be \emph{triangulable}
if $\fh$ is $\Ad$-conjugate to a triangular 
Lie subalgebra. 

2. A Lie subalgebra 
$\fh\subset\Lie(\Aut(X))$
is said to be \emph{locally integrable}
if it is filtered by an increasing sequence
of locally finite Lie subalgebras,
see \cite[Definitions~1.2]{AZ25b}. 
A solvable locally integrable 
Lie subalgebra $\fh\subset\Der(\kk[x,y])$
is triangulable, see 
\cite[Theorem~1.3]{AZ25b} 
 and the subsequent discussion. 
\edefis
\brem 
The Lie subalgebras $\fj_2^+$ and 
$\fj_2^-$ are conjugate by $\Ad_{\tau}$,
where $\tau\colon (x,y)\mapsto (y,x)$. Thus, 
a Lie subalgebra $\fh\subset\Lie(\Aut(\AA^2))$ 
is triangulable if and only if it is 
conjugate to a subalgebra of
\be\la{eq:j2} \fj_2^+=
 \kk[y]\frac{\p}{\p x}\oplus 
\kk x\frac{\p}{\p x}\oplus \kk y\frac{\p}{\p y}
\oplus\kk \frac{\p}{\p y}.\ee
\erem
Note that there exist 
solvable, non-triangulable Lie subalgebras of
$\Lie(\Aut(\AA^2))$,
see \cite[Remark~4.7]{AZ25b}. 
However, we have the following lemma. 
\blem \la{lem:crit-triangle}
 Let 
$\fh\subset\Lie(\Aut(\AA^2))$ be  
a weakly locally finite
solvable Lie subalgebra.
Then $\fh$ is triangulable and 
its derived length
is at most three.
\elem
\bproof Since $\fh$ is weakly locally finite 
and solvable, it
is contained 
in a locally integrable Borel subalgebra 
$\fb$ of $\Lie(\Aut(\AA^2))$,
see \cite[Lemma~3.2]{AZ25b}. 
By \cite[Theorem~3.4]{AZ25b}, 
$\fb=\Lie(B)$ for
some Borel subgroup $B$ of  
$\Aut(\AA^2)$. 
It is well known that any 
Borel subgroup of $\Aut(\AA^2)$
is conjugate to 
the triangular subgroup 
$\JONQ^+(\AA^2)$, 
see \cite[Theorem~1]{BEE16} for $\kk=\CC$ and 
\cite[Proposition~6.24]{CKRvS26} 
for an arbitrary algebraically 
closed field.
Therefore, $\fb$ is $\Ad$-conjugate 
to the subalgebra 
 of triangular derivations 
 $\fj_2^+=\Lie({\JONQ}^+(\AA^2))$.
 This proves the first statement. The second
 follows from the first; see Remark \ref{rem:fj-2}.1. 
\eproof
\brem\la{rem:Bass} 
In higher dimensions, an analogue 
of Lemma \ref{lem:crit-triangle}
is generally not valid. For example, 
by Bass \cite{Bas84},
the $\Ga$-subgroup 
$U=\exp(\kk\p)\subset\Aut(\AA^3)$, 
where $\p\in\Der(\kk[x,y,z])$ is
the locally nilpotent Nagata derivation, 
is not triangulable. Consequently, 
the abelian Lie subalgebra $\kk\p$
 is also not
triangulable. In fact, it is 
not even stably triangulable, 
 see \cite[Proposition~1.4.1]{KZ24}.
 This provides similar examples in all
 dimensions $n\ge 3$. 
\erem
\subsection{$\Lie(\Aut(\AA^2))$
 as a bigraded Lie algebra}
 Recall that the Lie algebra 
 $\Lie(\Aut(\AA^n))$ 
 coincides with the Lie algebra 
 $\Vec^{\mathrm c}(\AA^n)$ 
 of vector fields on $\AA^n$ 
 with constant divergence, 
 see, e.g., \cite[Proposition~15.7.2]{FK18}.
 The maximal ideal 
 $\Vec^0(\AA^n)\subset
 \Vec^{\mathrm c}(\AA^n)$ 
 of vector fields with zero divergence 
 coincides with the Lie algebra of 
 the normal subgroup 
 $\SAut(\AA^n)\subset\Aut(\AA^n)$
 of automorphisms whose 
 Jacobian equals $1$,
 see again \cite[Proposition~15.7.2]{FK18}. 
 
In what follows, 
we use the following notation,  
cf., e.g., \cite[Sec.~6.2]{KZ24} 
and \cite[(4)--(5)]{AZ25b}. 
\bnota
Let 
\[\Lambda=\{(a,b)\in\ZZ^2\,|
\,a,b\ge -1,\,\,\, (a,b)\neq (-1,-1)\}.\]
For $(a,b)\in\Lambda$, we define:
\be\la{eq:init} \p_{a,b}=
(b+1)x^{a+1}y^b\frac{\p}{\p x}
-(a+1)x^ay^{b+1}\frac{\p}{\p y}
\in{\Vec}^0(\AA^2)\ee
and 
\[{\rm bideg}(\p_{a,b})=(a,b).\]
This defines a bigrading on 
the Lie subalgebra 
$\Vec^0(\AA^2)$. 
Indeed, the following 
commutation relation
can be verified directly from the definition:
\be\la{eq:comm} 
[\p_{a,b},\p_{a',b'}]=
\det\left(\begin{matrix} 
a'+1&a+1\\b'+1&b+1\end{matrix}\right)
\p_{a+a', b+b'},\ee
where
\[(a,b),\,(a',b')\in\Lambda
\,\,\,\text{and} \,\,\, \p_{-1,-1}:=0.\]
All graded pieces of 
${\Vec}^0(\AA^2)$ 
are one-dimensional, 
and we have a decomposition
\[{\Vec}^0(\AA^2)=
\bigoplus_{(a,b)\in\Lambda} \kk\p_{a,b},
\quad\text{where}\quad
[ \kk\p_{a,b},  \kk\p_{c,d}]\subset  
\kk\p_{a+c,b+d},\] 
see \cite[Sec.~6.2]{KZ24}.
This bigrading respects 
the natural bigrading of $\kk[x,y]$, 
namely:
\be\la{eq:concord} 
\p_{i,j}(x^ky^l)=\det\left(\begin{matrix} 
k&i+1\\l&j+1\end{matrix}\right)
x^{k+i}y^{l+j}.\ee
We also have
\[{\Vec}^{\mathrm c}(\AA^2)
={\Vec}^0(\AA^2)
\oplus \kk\p_{\rm Eul}, 
\quad\text{where}\quad 
\p_{\rm Eul}=x\frac{\p}{\p x}+y\frac{\p}{\p y}\]
is the Euler derivation. Indeed, 
the Euler derivation spans the 
unique one-dimensional 
complement to ${\Vec}^0(\AA^2)$
in ${\Vec}^{\mathrm c}(\AA^2)$
because it has 
nonzero constant divergence. 
Therefore, every 
$\p\in{\Vec}^{\mathrm c}(\AA^2)$ 
admits a unique decomposition
\[\p=c_0\p_{\rm Eul}+
{\sum}_{i,j}c_{i,j} \p_{i,j},\]
where the sum is finite. 

The Lie algebra 
${\Vec}^{\mathrm c}(\AA^2)$ 
is bigraded with one-dimensional graded 
pieces, except
for the unique two-dimensional  
graded piece $\ft_2$ of weight $(0,0)$,
where
\be\la{eq:ss} \ft_2=\{\delta_{\alpha,\beta}\,|
\,(\alpha,\beta)\in\kk^2\}
\quad\text{with}\quad
\delta_{\alpha,\beta}=\alpha x\frac{\p}{\p x}
+\beta y\frac{\p}{\p y}\ee
is a maximal toral subalgebra
of ${\Vec}^{\mathrm c}(\AA^2)$.
The additional commutation relations 
are as follows:
\be\la{eq:comrel}
[\delta_{\alpha, \beta}, 
\delta_{\gamma, \eta}] 
= 0\quad\text{and}\quad
[\delta_{\alpha, \beta}, \p_{a,b}]= 
(\alpha a + \beta b) \p_{a, b}.\ee
Thus, every derivation $\p_{a, b}$ 
is an eigenvector of 
$\ad_{\delta_{\alpha, \beta}}\in
\End({\Vec}^0(\AA^2))$. 
\enota
\bdefi We call the lattice points 
$(-1,n)$ and $(m,-1)$ with $m,n\ge 0$
\emph{Demazure points}. 
They correspond to 
the homogeneous locally nilpotent
derivations
\[\p_{-1,n}=(n+1)y^n\frac{\p}{\p x}
\quad\text{resp.}\quad
\p_{m,-1}=-(m+1)x^m\frac{\p}{\p y},\]
and every homogeneous locally nilpotent
derivation $\p$ of $\kk[x,y]$ 
is proportional to one of them.
\edefi
\subsection{The technique of 
Newton polygons}
Here we introduce a piece of convex 
combinatorial geometry that will be 
used in several proofs. 
\bdefi\la{def:NP} Given a derivation 
$\p_0=\sum_{i,j} c_{i,j}\p_{i,j}
\in\Vec^{0}(\AA^2)$, 
the \emph{Newton polygon} 
$N(\p_0)$ is the
convex hull of 
\[\supp(\p_0):=\{(i,j)
\in\ZZ^2\,|\,c_{i,j}\neq 0\}.\]
 For a derivation $\p=\p_0+c_0\p_{\rm Eul}
\in\Vec^{\mathrm c}(\AA^2)$,
 where $c_0\neq 0$,
the Newton polygon $N(\p)$ is the
convex hull of
\[\supp(\p):=\supp(\p_0)\cup (0,0).\] 
\edefi
In the sequel, we use the 
following result.
\blem[{\rm \cite[Principle II]{vdE92}}]
\la{lem:vdV}
Let $\fh = \bigoplus_{i\in\ZZ} \fh_i$ 
be a $\ZZ$-graded Lie subalgebra of 
$\Lie(\Aut(X))$. If a nonzero derivation
$\p\in\fh$ is 
 locally finite of nonzero degree, 
 then the highest degree homogeneous 
 term of $\p$ with respect to this  grading
 is locally nilpotent. 
\elem
 \brem\la{rem:tot-ord}
The proof of this lemma remains 
valid verbatim 
if $\fh$ is graded by 
a totally ordered group, 
e.g. by a countable subgroup 
$\Gamma$ of $(\RR,+)$.
\erem
Using the grading on 
$\Vec^{\mathrm c}(\AA^2)$, 
we can deduce the following 
lemma.
\blem
\la{cor:loc-fin-hom-der} 
For every locally finite derivation
$\p\in \Vec^{\mathrm c}(\AA^2)$, 
all vertices of  the Newton polygon 
$N(\p)$ are among 
the lattice points 
$(-1,j), (i,-1)$ and $(0,0)$, 
where $i,j\ge 0$. 
If $\p$ is locally nilpotent, 
then all vertices 
of $N(\p)$ 
are Demazure points. 
\elem
\bproof  Let
$v=(i_0,j_0)$ be a vertex of $N(\p)$ 
distinct from $(0,0)$. There exists
a nonzero linear form 
$L=ax+by$ with 
integer coefficients such that
 $L$ attains its maximal value 
 $M$ on $N(\p)$  
 at the unique point $v$, 
 with $M>0$. 
 
Let $c\in\kk^*$. 
By setting $\deg_L(cx^ky^l)=ak+bl$, 
we obtain a $\ZZ$-grading on 
the polynomial ring $\kk[x,y]$,
and by setting $\deg_L(c\p_{i,j})=ai+bj$ 
and 
$\deg_L(\delta_{\xi,\eta})=0$ for all 
$\delta_{\xi,\eta}\in\ft_2$,
we define a $\ZZ$-grading 
on the Lie algebra 
 $\Vec^{\mathrm c}(\AA^2)$.
 
Suppose that 
$\p_{i,j}(x^ky^l)\neq 0$.
From \eqref{eq:concord} 
we deduce:
\be\la{eq:degrees}
{\deg}_L(\p_{i,j}(x^ky^l))=
{\deg}_L(x^ky^l)+{\deg}_L(\p_{i,j}).
\ee

Given $p\in \kk[x,y]$ (resp.
$\delta\in\Vec^{\mathrm c}(\AA^2)$),
we let $p_L$ (resp. $\delta_L$) 
denote the principal part
of $p$ (resp. of $\delta$) 
with respect to these gradings. 
According to our choice, 
$\p_L=c_0\p_{i_0,j_0}$, 
where $c_0\in\kk^*$.
Consequently,
$\deg_L(\p)=
\deg_L(c_0\p_{i_0,j_0})=M>0$. 
We claim that  
$\p_L$ is locally nilpotent. 
To prove the claim,
we follow the proof of 
Lemma \ref{lem:vdV} in
\cite{vdE92}.

Suppose the contrary. Then 
$\p_L^{(n)}(p)\neq 0$ for some 
$L$-homogeneous polynomial
$p=p_L\in\kk[x,y]$ 
and for all $n\ge 1$. Given our choice of $L$, 
we deduce from \eqref{eq:degrees}, 
by induction, that 
\be\la{eq:n}{\deg}_L(\p^{(n)}(p))
={\deg}_L(\p_L^{(n)}(p))=
nM+{\deg}_L(p),\ee
see formula $(*)$ in \cite[p.~863]{vdE92}.
Therefore, 
$(\p^{(n)}(p))_{n\ge 0}\subset\kk[x,y]$
is an infinite sequence of 
linearly independent polynomials. 
This contradicts the assumption
of local finiteness of $\p$, 
which proves the first assertion.

To prove the second assertion, suppose 
to the contrary that $\p$ 
is locally nilpotent
and that $v=(0,0)$ is a vertex of $N(\p)$. 
We can choose $L$ as before, 
with the sole difference
that $M=0$ is now the 
maximal value of $L|_{N(\p)}$
attained at the unique point:
the origin. 
Thus, the $L$-homogeneous leading part
$\p_L$ of $\p$ is a semisimple derivation
$\p_L=c\delta_{\alpha, \beta}$, 
where $c\in\kk^*$. 
Let us choose $k,l\in\NN$ such that 
\[\delta_{\alpha, \beta}(x^ky^l)=
(k\alpha+l\beta)x^ky^l\neq 0
\quad\text{and}\quad
{\deg}_L(x^ky^l)=ak+bl\neq 0.\]
It follows that, 
for all $n\ge 1$,
\[\p_L^{(n)}(x^ky^l)=
c^n(k\alpha+l\beta)^nx^ky^l\neq 0
\quad\text{and}\quad
{\deg}_L\left(\p_L^{(n)}(x^ky^l)\right)
={\deg}_L(x^ky^l)\neq 0.\]
Then $\p^{(n)}(x^ky^l)\neq 0$
as well for all $n\ge 1$,
see \eqref{eq:n}.
This is impossible, 
since we assume that $\p$ is 
locally nilpotent. 
\eproof
In what follows, we use the following 
combinatorial lemma.
\blem\la{lem:sum}
Let $\Pi_1$ and $\Pi_2$ be two 
convex polygons in $\RR^2$ and $L$ be
a linear form on $\RR^2$ such that
$L|_{\Pi_2}$ attains its maximum 
at a unique vertex $v$ of $\Pi_2$, 
and $L|_{\Pi_1}$ 
reaches its maximum either at 
a unique vertex 
$u_1$ of $\Pi_1$, 
or on a $1$-dimensional facet 
$F_1=[u_1,u_2]$ of $\Pi_1$. 
Consider the  convex hull $\Pi$ of 
$\Pi_1+\Pi_2$. Then $L|_{\Pi}$
attains its maximum at a unique vertex
$u_1+v$ of $\Pi$ in the first case, 
and on the $1$-dimensional facet 
$[u_1+v, u_2+v]$ of $\Pi$ 
in the second case.
\elem
\bproof Suppose we 
are in the second case; 
the argument in the first case
 is similar.  
Let  $u_1,\ldots,u_l$ and 
$v_1,\ldots,v_k$
be the vertices of $\Pi_1$ and 
$\Pi_2$, respectively, with $v_1=v$. 
Then $\Pi$ 
is the convex hull of the $u_i+v_j$. 
We have
\[\max(L|_{\Pi})={\max}_{i,j}\{L(u_i+v_j)\}
={\max}_i\{L(u_i)\} + {\max}_j\{L(v_j)\}=
L(u_1+v)=L(u_2+v).\] 
Since $L(u_i+v_j)<L(v+u_1)$ for $i>2$ 
or $j>1$, the lemma follows. 
\eproof
\brem If $L|_{\Pi_i}$ attains its maximum
on a $1$-dimensional facet $F_i$ 
for $i=1,2$, then $F_1$ and $F_2$ 
are parallel and $L|_{\Pi}$ 
attains its maximum
on a $1$-dimensional facet $F$ of $\Pi$
parallel to both of them. 
\erem
We apply this lemma to determine 
the Newton polygon of the bracket of two
derivations of $\kk[x,y]$.
\blem\la{cor:sum-newtons}
Let $\Pi_i=N(\p_i)$, $i=1,2$,
where $\p_1,\p_2\in\Der(\kk[x,y])$.
Suppose that $\Pi_1$ and $\Pi_2$
admit 
a common linear form $L$ 
satisfying the hypothesis of 
Lemma \ref{lem:sum}, where
$u_i=(k_i,l_i)$, $i=1,2$, 
 and $v=(m,n)$.
Suppose that, for some $i\in\{1,2\}$,  
we have
 \be\la{eq:ineqs} 
k_i+m \ge -1,\quad l_i+n \ge -1,\quad
(k_i+m,\,l_i+n)\neq (-1,-1),\ee
and
 \be\la{eq:non-vanish} 
[\p_{m,n}, \p_{k_i,l_i}]\neq 0, 
 \quad\text{i.e.}\quad 
 \det\left(\begin{matrix} 
k_i+1&m+1\\ l_i+1&n
+1\end{matrix}\right)
\neq 0\ee
(see \eqref{eq:comm}).
Then $u_i+v$ is a vertex of 
$N([\p_1,\p_2])$. 
If
\eqref{eq:ineqs}  
and \eqref{eq:non-vanish}
hold  for $i=1$ and $i=2$,
 then~$[u_1+v, u_2+v]$ is a
facet of 
$N([\p_1,\p_2])$. 
\elem
\bproof Assume again that we 
are in the second case of 
Lemma \ref{lem:sum}; 
the reasoning for the first case
 is analogous.  By Lemma \ref{lem:sum}, 
$u_1+v$ and $u_2+v$ 
are vertices of $\Pi=N(\p_1)+N(\p_2)$,
and $[u_1+v, u_2+v]$
 is a facet of $\Pi$. 
Suppose that \eqref{eq:ineqs}  and 
 \eqref{eq:non-vanish} are satisfied
for some $i\in\{1,2\}$. 
Since $N([\p_1,\p_2])\subset\Pi$,
 then 
 $u_i+v$ is a vertex of 
$N([\p_1,\p_2])$. If  \eqref{eq:ineqs}  and 
\eqref{eq:non-vanish}
hold for $i=1$ and $i=2$, then
 $[u_1+v, u_2+v]$ 
 is a facet of $N([\p_1,\p_2])$.
\eproof
\subsection{Commuting locally 
nilpotent derivations}
\la{ss:comm-lnd}
\blem \la{lem:centralizer-lnd}
Let $\p=p(y)\frac{\p}{\p x}$, where 
$p\in\kk[y]$ is nonzero, and let 
$\delta\in\Vec^{0}(\AA^2)$. 
Then $\p$ and
$\delta$
commute if and only if
one of the following holds:
\begin{itemize}
\item[$(i)$] $\deg(p)\ge 1$ and
$\delta=\kappa(y)\frac{\p}{\p x}$ 
for some $\kappa\in\kk[y]$;
\item[$(ii)$] $\deg(p)=0$ and
$\delta=\kappa(y)\frac{\p}{\p x}+
c_1\frac{\p}{\p y}$ 
for some $c_1\in\kk$ 
and some $\kappa\in\kk[y]$. 
\end{itemize}
Therefore,
$\Cent_{\Vec^{0}(\AA^2)}(\p)
\subset \fu^+_2
=\kk[y]\frac{\p}{\p x}\oplus \kk\frac{\p}{\p y}$ 
consists of locally nilpotent 
derivations. 
\elem
It is evident that if
either condition (i) or (ii) is satisfied, then
$\p$ and $\delta$ commute. For the converse,
we provide two different proofs:
the first combinatorial
and the second algebraic. 
\bproof[First proof] 
Note that the
Newton polygon $N(\p)$ 
is either a vertical segment $[v_0,v]$,
where  $v=(-1,d)$ and $v_0=(-1,d_0)$ with
$d_0<d=\deg(p)$,
or a singleton $\{v\}$.
The linear function $x$ restricted to 
$N(\delta)$
attains its maximal value either 
at a unique vertex, say $u_1$, or 
on a facet $[u_1,u_2]$ of $N(\delta)$, 
where for $u_1=(a,b)$ we choose 
the upper endpoint of this vertical segment. 
For  sufficiently small $\varepsilon>0$, 
the linear function 
$L=x+\varepsilon y$ reaches its maximal value 
$a+\varepsilon b$
on $N(\delta)$ at a unique vertex $u_1$. 
It attains
its maximal value $-1+\varepsilon d$
on $N(\p)$ at a unique vertex $v$.

If $[\p_{-1,d},\p_{a,b}]\neq 0$, 
then $a\ge 0$ and $L$
attains its maximal value on 
the Newton polygon $N([\p,\delta])$
at a unique vertex $u_1+v$ 
of $N([\p,\delta])$, 
see Lemma \ref{cor:sum-newtons}. 
However, if $\p$ and  $\delta$ commute, 
then $N([\p,\delta])$
is empty. This is  possible only if 
\[[\p_{-1,d},\p_{a,b}]=
(d+1)(a+1)\p_{a-1,b+d}= 0,\]
see \eqref{eq:comm}. 
The last equality
holds only if either $a=-1$ 
or $d=0$  and $(a,b)=(0,-1)$. Indeed,
if $d=0$ and $a\ge 0$,  then
the bracket vanishes only when
$\p_{a-1,b}=\p_{-1,-1}=0$.

Now, if $a=\max\{x|_{N(\delta)}\}=-1$, then 
$N(\delta)$ is contained in the vertical line 
$x+1=0$,
and therefore $\delta=q(y)\p/\p x$.
This corresponds to case (i).

Assume further that $a\ge 0$, $d=0$, 
and $u_1=(a,b)=(0,-1)$. 
In this case, we also have $u_2=(0,-1)$. 
Therefore, $\p=c_0\p/\p x$ and
$\delta=q(y)\p/\p x+c_1\p/\p y$
for some $q\in\kk[y]$, 
$c_0\in\kk^*$ and  $c_1\in\kk$. 
We are therefore in case (ii). 
In all cases, $\delta\in\fu^+_2$.
\eproof
\bproof[Second proof]
Let 
\[\delta=\alpha(x,y)\frac{\p}{\p x}
+\beta(x,y)\frac{\p}{\p y}, 
\quad\text{where}\quad
\frac{\p \alpha}{\p x}+\frac{\p\beta}{\p y}=0.\]
The equality $[\p,\delta]=0$ implies that 
\[p(y)\frac{\p\alpha}{\p x}-\beta p'(y)=0
\quad\text{and}\quad p(y) 
\frac{\p \beta}{\p x}=0.\]
Since $p\neq 0$, we obtain
\be\la{eq:first}
\beta=\beta(y)\in\kk[y]
\quad\text{and}\quad 
\alpha=\gamma(y)x+\kappa(y)
\quad\text{with}\quad \gamma,
\,\kappa\in\kk[y],\ee
where, according to the above, 
$\gamma$ and $\kappa$
satisfy the relations
\be\la{eq:second} p\gamma-\beta p'=0
\quad\text{and}\quad 
\gamma+\beta'=0.\ee
From the equalities \eqref{eq:second} 
we deduce that 
$(p\beta)'_y=
p\beta'+\beta p'=0$. 
Since $p,\beta\in k[y]$, 
the polynomial $p\beta$ is constant.
Thus, either $\beta=0$ or
$p$ and $\beta$ are both constants.

If $\beta=0$, then according 
to \eqref{eq:first} 
and \eqref{eq:second}, 
$\gamma=0$ and $\alpha=\kappa(y)$, 
 so that the pair
  \[(\p,\delta)
=\left(p(y)\frac{\p}{\p x}, \,
\kappa(y)\frac{\p}{\p x}\right)\]
satisfies (i) or (ii), depending on 
the value of $\deg(p)$. 

If $\beta=\text{const}\neq 0$, then 
$\gamma=0$,
$\deg(p)=0$, and thus the pair 
\[(\p,\delta)
=\left(c_0\frac{\p}{\p x},
 \kappa(y)\frac{\p}{\p x}
 +c_1\frac{\p}{\p y}\right)\]
satisfies (ii), where 
$c_0=p\in\kk^*$ 
and $c_1=\beta\in\kk$. 
It also follows 
that if $\beta=c_1\neq 0$, 
then $\deg(p)=0$, 
which completes the proof.
\eproof
\bcor\la{cor:centralizer-lnd} 
The centralizer  in $\Vec^{0}(\AA^2)$ 
of a nonzero locally 
nilpotent derivation 
consists of locally nilpotent derivations 
and is triangulable.
\ecor
\bproof Let 
$\p\in\Vec^{0}(\AA^2)\setminus\{0\}$ be 
locally nilpotent.
By Rentschler's theorem \cite{Ren68}, 
in a suitable coordinate system on $\AA^2$,
$\p$ is of the form $q(y)\frac{\p}{\p x}$, 
where $q\in\kk[y]\setminus\{0\}$. 
The result then follows
from Lemma  \ref{lem:centralizer-lnd}, 
since all $\delta$ 
as in items (i) and (ii) of this lemma are 
locally nilpotent triangular derivations.
\eproof
\brems $\,$
1. 
The assumption $\delta\in\Vec^0(\AA^2)$ 
in Corollary \ref{cor:centralizer-lnd} is essential:
the centralizer in $\Vec^{\mathrm c}(\AA^2)$ 
of a nonzero locally nilpotent derivation can 
contain derivations that are not locally nilpotent. 
Consider, for example,  the pair
$\p=\frac{\p}{\p x}$ and 
$\delta=\alpha\frac{\p}{\p x}+\beta y\frac{\p}{\p y}$
with $\alpha\in\kk$, $\beta\in\kk^*$. 

2. The conclusion of Corollary 
\ref{cor:centralizer-lnd} is not valid
in higher dimensions. For example, 
the semisimple derivation 
$x\frac{\p}{\p x}-y\frac{\p}{\p y}$
and the locally nilpotent derivation $\p/\p z$
commute and belong to $\Vec^{0}(\AA^3)$.
\erems
\bcor\la{cor:abel} 
Let $\fh\subset\fu^+_2$ 
be a commutative Lie subalgebra 
of dimension at least two
consisting of locally nilpotent derivations. 
Then either 
 \begin{itemize}
  \item[$(i)$]
  $\fh=\kk\frac{\p}{\p x}
  +\kk\frac{\p}{\p y}$,\\
or 
  \item[$(ii)$]
  $\fh\subset\kk[y]\frac{\p}{\p x}$.
   \end{itemize}
\ecor
\bproof Suppose that 
$\fh\not\subset\kk[y]\frac{\p}{\p x}$;
we will show that we are then in case 
(i).
Let us choose an arbitrary pair 
 $(\p_1,\p_2)$ of
linearly independent
derivations  in $\fh$ 
such that 
$E:=\span(\p_1,\p_2)\not
\subset\kk[y]\frac{\p}{\p x}$. Let us write
$\p_i=
p_i\frac{\p}{\p x}+c_i\frac{\p}{\p y}$,
where $p_i\in\kk[y]$ and $c_i\in\kk$
for $i=1,2$. By replacing one of these
derivations
with an appropriate linear combination, 
we may assume that $c_1=0$, and thus
$c_2\neq 0$. We claim that
$\deg(p_1)=0$, and consequently
$E=\kk\frac{\p}{\p x}
  +\kk\frac{\p}{\p y}$. 
Since every two-dimensional subspace of 
$\fh$ not contained in $\kk[y]\frac{\p}{\p x}$ 
equals $\kk\frac{\p}{\p x}\oplus\kk\frac{\p}{\p y}$,
 we are in case (i). 
 
To prove the claim, it suffices to note
that $[\p_2,\p_1]=c_2p_1'=0$. 
\eproof
\subsection{Nilpotent Lie subalgebras 
generated by locally 
nilpotent derivations}\la{ss:nilp}
Note that any locally nilpotent derivation 
has divergence zero. 
Therefore, a Lie subalgebra
of $\Vec^{\mathrm c}(X)$ 
generated by locally 
nilpotent derivations is contained in 
$\Vec^{0}(\AA^2)$.
In this subsection we prove
the following proposition.
\bprop\la{thm:center} 
Suppose that a solvable Lie subalgebra 
$\fh\subset \Vec^{0}(\AA^2)$ contains a
nonzero locally nilpotent 
derivation $\p_0$ and
a nonzero central element $\eta_0$.
Then $\fh$ is triangulable
 and consists 
of locally nilpotent derivations. 
If $\fh$ is finitely generated, then
$\fh$ is locally finite. 
\eprop
\bproof
Since $[\p_0,\eta_0]=0$ and $\p_0$ is 
locally nilpotent, $\eta_0$ 
is also locally nilpotent
by Corollary \ref{cor:centralizer-lnd}.
By Rentschler's theorem \cite{Ren68}, 
in a suitable coordinate 
system on $\AA^2$,
$\eta_0$ is of the form 
$\eta_0=q(y)\frac{\p}{\p x}$, 
where $q\in\kk[y]$.  
Since every $\p\in\fh$  
commutes with $\eta_0$, according to
Lemma \ref{lem:centralizer-lnd},
$\p\in\fu^+_2$ 
is locally nilpotent.
Therefore, $\fh\subset\fu^+_{2}$. 
Since every element of $\fu_2^+$ 
is locally nilpotent, $\fh$ consists 
of locally nilpotent derivations.
If $\fh$ is finitely generated, then
$\fh\subset\fu^+_{2,d}$  
for some $d>0$. 
In the latter case, $\fh$  is locally finite,
since $\fu^+_{2,d}$ is, see \eqref{eq:fu-2}.
\eproof
The following corollary is immediate.
Indeed, every nilpotent Lie algebra 
has nontrivial center.
\bcor\la{cor:nilp}
A nilpotent Lie subalgebra 
$\fh\subset \Vec^{0}(\AA^2)$ 
that contains a
nonzero locally nilpotent 
derivation  is triangulable and consists 
of locally nilpotent derivations. 
\ecor
\subsection{Solvable Lie subalgebras 
generated by locally 
nilpotent derivations}\la{ss:solv}
In this subsection we 
prove the following theorem.
\bthm\la{thm:unip}
Let $\fh$
be a solvable Lie subalgebra 
of $\Vec^{0}(\AA^2)$ 
generated by locally 
nilpotent derivations.
Then $\fh$ is triangulable. 
If $\fh$
is generated by 
a finite set of locally 
nilpotent derivations, 
then $\fh=\Lie(U)$ 
for a unipotent algebraic
 subgroup
$U\subset\Aut(\AA^2)$, 
and therefore $\fh$ is locally finite. 
In general, $\fh=\Lie(U)$, 
where $U$
is a nested unipotent ind-subgroup 
of $\Aut(X)$. Consequently, $\fh$
consists of locally 
nilpotent derivations derivations 
and is weakly locally finite.
\ethm
Before passing to the proof, 
we deduce the following corollary.
\bcor\la{cor:lnd}
Let $\fh\subset\Vec^{\mathrm c}(\AA^2)$
be a $J$-saturated solvable Lie subalgebra 
generated by 
locally finite derivations. 
Consider the Lie subalgebra $\fh_{\rm lnd}$
of $\fh$ generated by all locally nilpotent 
derivations from $\fh$. Then $\fh_{\rm lnd}$
is an ideal of $\fh$ consisting 
of locally nilpotent derivations.
\ecor
\bproof According to 
Lemma \ref{lem:ideal},
$\fh_{\rm lnd}$
is an ideal of $\fh$.
By virtue of
Theorem  \ref{thm:unip},  $\fh_{\rm lnd}$
consists of locally nilpotent derivations.
\eproof
The proof of Theorem  \ref{thm:unip}
is preceded by Lemmas 
\ref{lem:lnds}--\ref{cor:primitive}.
\blem\la{lem:lnds} Consider
a solvable Lie subalgebra $\fh$ of 
$\Vec^{0}(\AA^2)$  that
contains a nonzero locally 
nilpotent derivation $\p$.
Suppose that $\fh^{(d-1)}\neq 0$ 
and $\fh^{(d)}= 0$, 
where $d\ge 1$. Then  $\fh^{(d-1)}$ 
consists of locally nilpotent derivations
and is triangulable.
\elem
\bproof It suffices to show 
that $\fh^{(d-1)}$ 
contains a nonzero locally 
nilpotent derivation, 
say $\zeta$. 
Indeed, since $\fh^{(d-1)}$ is commutative,
$\zeta$ is central in $\fh^{(d-1)}$, 
and then the conclusion follows
from Corollary \ref{cor:centralizer-lnd}. 

To show the existence of $\zeta$, 
we proceed by induction on $d$. 
If $d=1$, then setting $\zeta=\p$ 
we are done.   
If $d>1$, suppose by induction that 
$\fh^{(j-1)}$,
where $1\le j\le d-1$, contains 
a nonzero locally nilpotent derivation
$\p_j$. 
Since $\fh^{(j-1)}$ is not abelian, there is 
a derivation $\delta\in\fh^{(j-1)}$ such that
$\p_j$ and $\delta$ do not commute.
Let $m>1$ be such that 
$\eta:=\ad^{(m-1)}_{\p_j}(\delta)\neq 0$
and 
$\ad_{\p_j}(\eta)=
\ad^{(m)}_{\p_j}(\delta)=0$.
Thus, $\eta\in \fh^{(j)}$ is nonzero 
and commutes with $\p_j$. 
By Corollary \ref{cor:centralizer-lnd},
$\eta$ is locally nilpotent. 
This completes the induction. 
Thus, $\fh^{(d-1)}$ 
contains a nonzero locally 
nilpotent derivation.
\eproof
\bnota\la{not:temporary}
In what follows, we denote
by $\fh$ a solvable Lie subalgebra 
of $\Vec^{0}(\AA^2)$ 
with derived length $d\ge 1$,
generated by the locally 
nilpotent derivations $\p_i$,
$i=1,2,\ldots$. 
 \enota
We use the following terminology.
\bdefi
A derivation 
$\p=a(x,y)\frac{\p}{\p x}+b(x,y)
\frac{\p}{\p y}\in\Der(\kk[x,y])$
 is called \emph{irreducible}
if $a$ and $b$ are coprime, 
and \emph{reducible} otherwise.
\edefi
Recall the following lemma,
see e.g. \cite[Corollary~2.13]{Fre17}.
\blem\la{lem:Fre}
A locally nilpotent derivation 
$\p\in\Der(\kk[x,y])$ is irreducible
(resp., reducible) 
if and only if,
in appropriate coordinates 
$(u,v)$ of $\AA^2$,  
we have $\p=\p/\p v$
(resp.,
$\p=p(u)\p/\p v$, where
$p\in\kk[u]$ and $\deg(p)\ge 1$).
\elem
\blem\la{cor:triangle} 
Suppose that $\fh^{(d-1)}$ contains 
a reducible derivation $\zeta$. 
Then  $\fh$ is triangulable. 
More precisely, if $\zeta=p(u)\p/\p v$
in appropriate coordinates 
$(u,v)$ on $\AA^2$, where
$p\in\kk[u]$ and $\deg(p)\ge 1$, then
$\fh \subset \fu^+_2$ 
(with respect to these coordinates).
\elem
\bproof 
Suppose that 
$\zeta=p(u)\p/\p v\in\fu^+_2$ 
in appropriate coordinates 
$(u,v)$ on $\AA^2$, 
where $\deg(p)>0$. 
Then, by Lemma 
\ref{lem:centralizer-lnd} (case (i)), 
$\fh^{(d-1)}\subset\kk[u]\p/\p v$. 

For every locally nilpotent generator 
$\p_i$ of $\fh$,
the restriction of $\ad_{\p_i}$  
to the ideal $\fh^{(d-1)}$ is
a locally nilpotent derivation 
of $\fh^{(d-1)}$,
whose kernel is nonzero. Let 
$\eta_i\in\ker(\ad_{\p_i}|_{\fh^{(d-1)}})$
be nonzero, where, by the preceding,
$\eta_i=q_i(u)\p/\p v$ 
with a nonzero $q_i\in\kk[u]$. 
Since $\p_i$ and $\eta_i$ commute,
with respect to the coordinates $(u,v)$
we have
$\p_i\in \fu^+_2$ for all $i=1,2,\ldots$, 
see again Lemma 
\ref{lem:centralizer-lnd}.
Therefore, $\fh\subset  \fu^+_2$.
\eproof
\blem\la{lem:primitive}
Suppose that all nonzero derivations 
of $\fh^{(d-1)}$
are irreducible. Then $\fh^{(d-1)}$ 
is $\Ad$-conjugate
to a subalgebra of the Lie algebra
$\kk\frac{\p}{\p u}\oplus\kk\frac{\p}{\p v}$,
where $(u,v)$ is a suitable 
coordinate system on $\AA^2$.
\elem
\bproof By Lemma \ref{lem:lnds}, 
in suitable coordinates $(u,v)$ of $\AA^2$
we have
\[\fh^{(d-1)}\subset \fu^+_2=
\kk[u]\frac{\p}{\p v}\oplus\kk\frac{\p}{\p u}.\] 
Suppose that
$\fh^{(d-1)}\subset\kk[u]\frac{\p}{\p v}$.
Then, in fact, $\fh^{(d-1)}=\kk\frac{\p}{\p v}$,
since it contains
no reducible derivation. 
According to Rentschler's theorem,
this is the case if $\dim \fh^{(d-1)}=1$.

Furthermore, suppose that 
$\dim \fh^{(d-1)} \ge 2$ and 
$\fh^{(d-1)}\not\subset\kk[u]\frac{\p}{\p v}$. 
Then, by Corollary \ref{cor:abel},
\[\fh^{(d-1)}=
\kk\frac{\p}{\p u}\oplus\kk\frac{\p}{\p v}.\] 
\eproof
\blem\la{cor:primitive}
Under the assumptions of 
Lemma \ref{lem:primitive},
$\fh$ is triangulable. 
\elem
\bproof 
Let $\widetilde \fh = 
\{\ad_{\p}|_{\fh^{(d-1)}} \,|\,\p \in \fh\}$.
Then  $\widetilde\fh$ is 
a Lie subalgebra of $\End(\fh^{(d-1)})$,
where $\dim(\fh^{(d-1)})\le 2$.
Since $\fh$ is solvable, 
so is $\widetilde\fh$.
So, Lie's Theorem implies 
that $\widetilde\fh$ is 
a triangulable Lie subalgebra of 
$\End( \fh^{(d-1)})$.
Since $\widetilde\fh$ is triangulable 
and is generated 
(as a Lie algebra) 
by the nilpotent endomorphisms 
$\ad_{\p_i}|_{\fh^{(d-1)}}$,
it follows that $\widetilde\fh$ consists of
nilpotent endomorphisms.
Then Engel's Theorem gives a common 
nonzero kernel element of $\widetilde\fh$, 
so the center of $\fh$ is nonzero.
The conclusion then follows 
from Proposition \ref{thm:center}. 
\eproof
\bproof[Proof of Theorem \ref{thm:unip}]
By Lemmas \ref{cor:triangle}--\ref{cor:primitive}, 
$\fh$ is triangulable. 
Therefore, we can assume that 
$\fh\subset \fu^+_2$. If $\fh$ 
is generated by a finite set
of locally nilpotent derivations $\p_i$, 
then
$\p_i\in\fu^+_{2,d}$ for all $i$ and 
a suitable $d\gg 1$. 
Consequently, $\fh\subset\fu^+_{2,d}$ is 
a finite-dimensional nilpotent subalgebra 
consisting of locally nilpotent derivations.
Thus, $U:=\exp(\fh)\subset\exp(\fu^+_{2,d})$
is a unipotent algebraic subgroup 
with Lie algebra $\Lie(U)=\fh$. 
In particular, $\fh$ is locally finite. 
In general, $\fh=\Lie(U)$, where $U=\exp(\fh)$
is a nested unipotent subgroup of 
$U_2^+=\exp(\fu_2^+)$. 
\eproof
\subsection{Spectral decomposition 
related to
a semisimple derivation} 
We begin with the following observations. 
\brem\la{rem:ss}
Let $s\in\Vec^{\mathrm c}(\AA^2)$ 
be a semisimple derivation.
The toral subalgebra 
$\kk s\subset\Vec^{\mathrm c}(\AA^2)$
is contained in $\Lie(T_{\min}(s))$
for a smallest algebraic torus 
$T_{\min}(s)\subset\Aut(\AA^2)$,
see Lemma~\ref{lem:toral}.
If $s$ is not algebraic, 
we have $\dim(T_{\min}(s))=2$.
In this case, according to 
Bia\l{}ynicki-Birula's theorem 
\cite[I]{Bia66-67}, 
$T_{\min}(s)$ is conjugate 
in $\Aut(\AA^2)$
to the 2-torus $\TT\subset{\rm GL}(2,\kk)$
consisting of diagonal matrices, 
and $\kk s$ is conjugate
to a Lie subalgebra 
$\kk\delta_{\alpha,\beta}
\subset\ft_2=\Lie(\TT)$,
where $\alpha,\beta\in\kk$ are 
nonzero with 
$\alpha/\beta\in\kk\setminus\QQ$. 

In the case where $T_{\min}(s)$
is a one-dimensional 
algebraic torus,
that is, $s$ is algebraic, 
$T_{\min}(s)$ is conjugate to a 
subtorus of the diagonal 
$2$-torus $\TT$, see 
\cite[II]{Bia66-67} or
\cite{Gut62}.
Therefore, $\kk s$ is
$\Ad$-conjugate
to $\kk\delta_{m,n}$ for
a pair $(m,n)$ of 
coprime integers.
\erem

Recall the following lemma; 
see, for example,
 \cite[Lemma~4.8]{AZ25b}.
\blem
\la{lem-simple}
Let $V$ be a vector space  over $\kk$ 
and $A$
a diagonalizable endomorphism of $V$.
Consider the decomposition 
$V = \bigoplus_\lambda V_\lambda$, 
where the $V_\lambda$ 
are the eigenspaces 
of $A$ associated with 
the different eigenvalues 
$\lambda\in\kk$. 
Let a subspace $U \subset V$ 
be invariant under $A$, and let 
$u=\sum_\lambda u_\lambda \in U$, 
where $u_\lambda \in V_\lambda$. 
Then $u_\lambda \in U$ 
for all $\lambda$.
\elem
Given a pair $(\alpha,\beta)\in
\kk^2\setminus \{(0,0)\}$,
the action of the semisimple 
derivation 
$\ad_{\delta_{\alpha,\beta}}$ on 
${\Vec}^{\mathrm c}(\AA^2)$ 
leads to
a spectral decomposition
\[{\Vec}^{\mathrm c}(\AA^2)=
\bigoplus_{\lambda\in\kk} 
E_\lambda(\delta_{\alpha, \beta}),\]
where 
$E_\lambda(\delta_{\alpha, \beta})$ 
is the eigenspace
of $\ad_{\delta_{\alpha, \beta}}$ 
associated with the eigenvalue 
$\lambda$. 
The eigenspace 
$E_0(\delta_{\alpha, \beta})$ 
coincides 
with the centralizer of 
$\delta_{\alpha, \beta}$ in 
${\Vec}^{\mathrm c}(\AA^2)$. 
According to Lemma \ref{lem-simple},
for a vector subspace 
$V\subset{\Vec}^{\mathrm c}(\AA^2)$
 invariant under 
$\ad_{\delta_{\alpha,\beta}}$,
we have
\be\la{eq:dir}
V=\bigoplus_{\lambda\in\kk}\left(V\cap 
E_\lambda(\delta_{\alpha,\beta})\right),
 \ee
 where the sum is finite if $V$ 
 is finite-dimensional.
In particular, if $\delta_{\alpha,\beta}\in\fh$, 
where $\fh$ is a Lie subalgebra of 
${\Vec}^{\mathrm c}(\AA^2)$,
then we have 
\be\la{eq:dir-dec}
\fh=\bigoplus_{\lambda\in\kk}\left(\fh\cap 
E_\lambda(\delta_{\alpha,\beta})\right).
 \ee
 Note that any $\p\in\fh$ has only
 a finite number 
 of nonzero homogeneous components, 
 and all of them are eigenvectors of 
 $\ad_{\delta_{\alpha,\beta}}$. 
 Consequently,
the spectral decomposition
$\p=\sum_\lambda \p_\lambda$
associated with 
$\delta_{\alpha,\beta}$
is also finite. Moreover,
$\p_{\lambda}\in \fh\cap
E_\lambda(\delta_{\alpha,\beta})$ 
is the sum of the homogeneous components
of $\p$ that are eigenvectors of 
 $\ad_{\delta_{\alpha,\beta}}$ 
corresponding to the eigenvalue
  $\lambda$. 
\subsection{Solvable Lie subalgebras 
generated by locally finite
derivations}\la{ss:3.8}
In this subsection, we prove 
Theorem \ref{main-thm}. 
Since the triangular Lie subalgebra 
$\fj^+_2\subset\Vec^{\mathrm c}(\AA^2)$
is weakly locally finite --- meaning that
every finitely generated subalgebra of 
$\fj^+_2$ is locally finite --- it suffices to
prove the following result. 
\bthm
\la{thm:rk-1} 
Every solvable Lie subalgebra 
$\fh\subset\Vec^{\mathrm c}(\AA^2)$ 
generated by locally finite
derivations is triangulable.
\ethm
In the sequel, we adopt 
the following convention. 
\bconv\la{conv} 
By Proposition \ref{prop:max}
we may assume, and shall assume,
 that $\fh$ in Theorem \ref{thm:rk-1} is
 maximal 
among solvable subalgebras of 
$\Vec^{\mathrm c}(\AA^2)$ 
generated 
by locally finite derivations 
and is, therefore,
 $J$-saturated
 and generated by algebraic 
 Lie subalgebras $\fh_k$,  $k=1,2,\ldots$
 (see Proposition \ref{prop:alg-toral}). 
 By decomposing each $\fh_k$ into
 a semidirect product of 
 a maximal nilpotent ideal $\fm_k$
 and an algebraic torus $\ft_k$ of 
 rank at most two, 
 we deduce that $\fh$ is generated 
 by two sequences, the first consisting 
 of locally nilpotent  derivations $\p_i$
and the second consisting 
 of semisimple algebraic 
 derivations $\delta_j$. 
 
Let $\fh_{\rm lnd}$ 
 be the Lie subalgebra of $\fh$ 
 generated by all locally nilpotent  
 derivations of $\fh$. By Theorem
 \ref{thm:unip}, $\fh_{\rm lnd}$ 
 is triangulable
 and consists of locally 
 nilpotent derivations.
 We therefore further assume
that $\fh_{\rm lnd}\subset\fu^+_2$.
 \econv
Since $\fh$, as defined in 
Convention \ref{conv},
is generated by $\fh_{\rm lnd}\subset \fu^+_2$ 
and the semisimple derivations 
from $\fh$, to prove Theorem  \ref{thm:rk-1}
it suffices to show that $\delta\in\fj^+_2$
for any semisimple derivation $\delta\in\fh$. 
Indeed, we then also have 
$\fh\subset\fj^+_2$, and thus
$\fh$ is triangular.  In fact, the 
inclusion $\delta\in\fj^+_2$, that 
holds in most cases, 
in certain situations
does not necessarily hold; 
we will  then adopt
a different approach. 
\blem\la{lem:invar} Let $\delta\in\fh$ be 
a semisimple derivation. 
Then $\ad_{\delta}$ 
leaves invariant every 
derived ideal 
$(\fh_{\rm lnd})^{(k)}$ 
of $\fh_{\rm lnd}$. 
\elem
\bproof
Since $\fh$ 
is assumed to be maximal 
among solvable subalgebras of 
$\Vec^{\mathrm c}(\AA^2)$ 
generated 
by locally finite derivations, 
$\fh$ contains 
$\Lie(T_{\min}(\delta))$,
see Proposition \ref{prop:alg-toral}
and its proof. 
By  Lemma \ref{lem:toral},
$\Ad_{T_{\min}(\delta)}$ leaves
invariant the ideal $\fh_{\rm lnd}$ 
of $\fh$.
Since $\Ad_{T_{\min}(\delta)}$
acts on $\fh_{\rm lnd}$ 
by automorphisms of this Lie algebra,
it preserves every derived ideal 
$(\fh_{\rm lnd})^{(k)}$ 
of $\fh_{\rm lnd}$. Since
$\delta\in\Lie(T_{\min}(\delta))$,
the endomorphism $\ad_\delta$
also preserves these ideals. 
\eproof
\brem\la{rem:rent} We have 
$\fu^+_2=\kk[y]\frac{\p}{\p x}
+\kk\frac{\p}{\p y}$ 
and
$[\fu^+_2,\fu^+_2]=\kk[y]\frac{\p}{\p x}$.
 It follows 
that either
$ \fh_{\rm lnd}\subset \fu^+_2$ is abelian, 
or otherwise 
the derived ideal $(\fh_{\rm lnd})^{(1)}$ 
is nonzero and 
contained in
$\kk[y]\frac{\p}{\p x}$. 
\erem
\blem\la{lem:spec}
Consider a locally nilpotent derivation 
$\p=p(y)\frac{\p}{\p x}$, 
where $p\in\kk[y]$ is nonzero. 
Let $\delta$ be 
a semisimple derivation of $\kk[x,y]$
such that $[\delta,\p]=\lambda\p$. 
Then $\delta\in\fj^+_2$.
\elem
\bproof Set $d=\deg(p)$.
The Newton polygon
$N(\p)$ is either a vertical segment
whose upper vertex is $v_1=(-1,d)$, 
 or the singleton $v_1=(-1,d)$.
For a  sufficiently small $\varepsilon>0$,
the linear function 
$L_\varepsilon=x+\varepsilon y$
attains its maximal value on $N(\p)$ at
a unique vertex $v_1$.

The function $x|_{N(\delta)}$ 
attains its maximal value either 
at a unique vertex $u_1=(k,l)$ 
or at a vertical facet $[u_1,u_2]$. 
In the latter case, 
we denote by $u_1$ the upper endpoint
of this facet. In both cases, 
for a  sufficiently small $\varepsilon>0$,
the linear function 
$L_\varepsilon |_{N(\delta)}$ 
attains its maximal value 
at a unique vertex $u_1$. 

Note that $k\ge 0$.
Indeed, if $k=-1$, then
$N(\delta)$ is contained in 
the vertical line $x+1=0$,
and therefore
$\delta\in\kk[y]\frac{\p}{\p x}$ 
is locally nilpotent, 
contrary to our choice. 

Since $\delta$ is locally finite, 
we deduce from
Lemma \ref{cor:loc-fin-hom-der} that
$u_1\in\{(0,0),\,(k,-1)\}$. If $k=0$, then, 
again by 
Lemma \ref{cor:loc-fin-hom-der}, 
$N(\delta)$ is contained 
in a certain convex 
quadrilateral $\Pi_n$ 
with vertices 
\[(-1,n),\,(-1,-1),\,(0,-1),\,(0,0),\]
where $n\ge 0$.
It follows that $\delta\in\fj^+_2$
(see \eqref{eq:j2}), as desired. 

Let us show that $k=0$. 
Indeed, if $k>0$, then $u_1=(k,-1)$.
By \eqref{eq:comm},
we obtain:
\be\la{eq:crochets}
[\p_{-1,d},\p_{k,-1}]=
(d+1)(k+1)\p_{k-1,d-1}
\neq 0.
\ee
According to 
Lemma \ref{cor:sum-newtons},
$L_\varepsilon |_{N([\p,\delta])}$ 
attains its maximal value 
at a unique vertex 
$u_1+v_1$ of $N([\p,\delta])$. 
Since the leading homogeneous 
components already have 
a nonzero bracket 
\eqref{eq:crochets}, 
the Newton polygon of 
$[\p,\delta]$ is nonempty.
 Hence $[\delta,\p]\neq0$, 
 so $\lambda\neq0$.
By our assumption,
$[\delta,\p]=\lambda\p$. 
Since $\lambda\neq 0$, 
we have $N([\p,\delta])=N(\p)$.
 It follows that $u_1+v_1=v_1$, 
 that is, $u_1=(0,0)$. 
 This contradicts 
 the equality $u_1=(k,-1)$.
Thus, $k=0$, and consequently
$\delta\in\fj^+_2$.
\eproof
\blem\la{lem:2-cases}
Suppose that 
$\fh_{\rm{lnd}}\subset\fu^+_2$ 
is either noncommutative, 
or nonzero 
and contained in 
$\kk[y]\frac{\p}{\p x}$.
Then $\fh\subset\fj^+_2$.
\elem
\bproof
In the first case, 
$(\fh_{\rm{lnd}})^{(1)}\subset
 \kk[y]\frac{\p}{\p x}$
is nonzero (see Remark \ref{rem:rent}), 
and in the second case, 
$\fh_{\rm{lnd}}
\subset\kk[y]\frac{\p}{\p x}$
is nonzero. 
For any semisimple derivation 
$\delta\in \fh$, 
the endomorphism $\ad_\delta$
preserves both $\fh_{\rm{lnd}}$ and 
$(\fh_{\rm{lnd}})^{(1)}$, 
see Lemma \ref{lem:invar}.
For the restriction 
$\ad_{\delta}|_{\fh_{\rm{lnd}}}$
(resp. 
$\ad_{\delta}|_{(\fh_{\rm{lnd}})^{(1)}}$),
there exists
$\lambda\in\Spec(\ad_{\delta}|_{\fh})$ 
and a nonzero derivation 
$\p_\lambda\in\fh_{\rm{lnd}}$
(resp. $\p_\lambda\in (\fh_{\rm{lnd}})^{(1)}$)
such that 
$\p_\lambda=p(y)\frac{\p}{\p x}
\in\kk[y]\frac{\p}{\p x}$
and
$[\delta,\p_\lambda]=\lambda\p_\lambda$. 
The conclusion then follows 
from Lemma \ref{lem:spec}.
\eproof
\brem\la{rem:abel}
Suppose that neither of the hypotheses
of Lemma \ref{lem:2-cases} 
holds for
$\fh_{\rm{lnd}}\subset\fu^+_2$, 
that is, suppose it is
commutative and not 
contained in $\kk[y]\frac{\p}{\p x}$.
Then $\fh_{\rm{lnd}}\subset
\kk\frac{\p}{\p x}+\kk\frac{\p}{\p y}$,
see  Corollary \ref{cor:abel}.
\erem
\blem\la{lem:aff}
Suppose that 
$\fh_{\rm{lnd}}=
\kk\frac{\p}{\p x}+\kk\frac{\p}{\p y}$.
Then $\fh$ is locally finite and triangulable. 
\elem
\bproof
Let us determine all derivations
\[\delta=p\frac{\p}{\p x}+q\frac{\p}{\p y}
\quad\text{with}\quad p,\,q\in\kk[x,y]\]
such that $\ad_\delta$ preserves
$\fh_{\rm{lnd}}$.
We have
\[
\Big[\delta,\frac{\p}{\p x}\Big]=
-\frac{\p p}{\p x}\frac{\p}{\p x}
-\frac{\p q}{\p x}\frac{\p}{\p y}
\quad\text{and}\quad
\Big[\delta,\frac{\p}{\p y}\Big]=
-\frac{\p p}{\p y}\frac{\p}{\p x}
-\frac{\p q}{\p y}\frac{\p}{\p y}.
\]
Since these two brackets belong to $\fh_{\rm{lnd}}=
\kk\frac{\p}{\p x}+\kk\frac{\p}{\p y}$, 
each coefficient must be constant. Thus,
$p$ and $q$ are affine linear, and therefore
\[\delta=
(ax+by+c)\frac{\p}{\p x}
+(dx+ey+f)\frac{\p}{\p y}\in{\rm aff}_2
:=\Lie(\Aff(\AA^2)).
\]
Since $\ad_{\fh}$ preserves the ideal 
$\fh_{\rm{lnd}}$ of $\fh$, 
we conclude that $\fh\subset {\rm aff}_2$ 
is locally finite. As it is solvable, 
$\fh$ is contained in a Borel subalgebra of 
${\rm aff}_2$. Every Borel subalgebra of 
${\rm aff}_2$ is conjugate to 
${\rm aff}^+_2:={\rm aff}_2\cap\fj^+_2$;
consequently, $\fh$ is triangulable.
\eproof
\blem\la{lem:lin} Suppose that  
$\fh_{\rm{lnd}}=\kk\p_0$,
where $\p_0=a_0\frac{\p}{\p x}+
b_0\frac{\p}{\p y}\neq 0$
with $a_0, b_0\in\kk$. 
Then $\fh$ is  triangulable.
\elem
\bproof
Since $\ad_{\fh}$
preserves $\fh_{\rm{lnd}}=\kk\p_0$, 
for any semisimple derivation
$\delta\in\fh$ we have 
$[\delta, \p_0]=\lambda\p_0$, 
where $\lambda\in\Spec(\ad_\delta)$. 
After an appropriate
 linear change of variables 
$(x,y)\mapsto (u,v)$, we may assume that 
$\p_0=\frac{\p}{\p u}$. 
In these new variables,
we have $\delta\in\fj^+_2$, 
see Lemma \ref{lem:spec}. 
Since $\fh$ is generated by $\p_0$ 
and by semisimple derivations, 
it follows
that $\fh\subset\fj^+_2$ 
(in the new coordinates). 
The lemma follows. 
\eproof
\bcor\la{cor:1-lnd} 
If $\fh$ contains at least one nonzero 
locally nilpotent derivation, then $\fh$ 
is triangulable. 
\ecor
\bproof In view of Remark \ref{rem:rent},
this follows immediately from 
Lemmas \ref{lem:2-cases},
\ref{lem:aff}, and \ref{lem:lin}. 
\eproof
 The following lemma is inspired by
\cite[Corollary~3.4]{FZ05}. It is 
adapted to our  more general setting of
Lie algebras of derivations.
\blem\la{lem:1-lnd} 
Let $\fh\subset\Lie(\Aut(X))$ 
be a Lie subalgebra.
Suppose that $\fh$ contains a pair
of derivations 
$\delta$ and $\p$ that do not commute,
where $\delta$ is semisimple and 
$\p$ is locally finite. 
Then $\fh$ contains a nonzero 
locally nilpotent derivation. 
\elem
\bproof 
Since $\delta$ and $\ad_\delta$ 
are locally finite, 
there exists a finite-dimensional
vector subspace $V\subset\fh$ 
 containing $\p$
and stable under the action 
of $\ad_\delta$.  Since $\delta$
is semisimple, $V$ possesses
a basis consisting of
eigenvectors of $\delta$, 
see \eqref{eq:dir}. 
Let $\sigma=\Spec(\ad_\delta|_V)
\subset\kk$, and let
$\p=\sum_{\lambda\in\sigma} \p_\lambda$
be the spectral decomposition of $\p$ 
with respect to $\ad_\delta|_V$,
where $\p_\lambda\in V$ 
and 
$[\delta,\p_\lambda]=\lambda \p_\lambda$.
Since $\p$ and $\delta$ do not commute, 
the finite set $\sigma\subset \kk$
 is not reduced to the singleton 
$\{0\}$.

Let $\QQ\subset K$ be 
the finite extension  
of $\QQ$ in $\kk$ spanned 
by $\sigma$.
The number field $K$
can be embedded
in $\CC$ as a subfield. 
Let $\cN$ be the convex hull  of 
$\sigma$ in $\CC$.
The vertices of $\cN$
belong to $\sigma$. 
Let us choose a linear form
$L\colon\CC=\RR^2\to\RR$ 
such that
$L|_\sigma$ does not 
vanish identically and, up to 
a sign change, 
attains its positive 
maximal value on $\cN$
at a unique vertex 
$\lambda_{\max}\in\sigma$ 
of $\cN$. 

Consider the finitely generated Lie subalgebra 
$\fl_V=\langle \p_\lambda\,|
\,\lambda\in \sigma\rangle_{\Lie}$. 
Since every $\p_\lambda$ belongs to
$\fh$, the Lie algebra
$\fl_V$ is a subalgebra of $\fh$.
It is invariant under $\ad_\delta$,
and $\Spec(\delta|_{\fl_V})$ is contained in
the $\ZZ$-submodule $\Lambda$ of $\CC$ 
generated by $\sigma$. 
 The image $\Gamma:=L(\Lambda)$ 
 is a finitely generated subgroup 
 of $(\RR,+)$. It is therefore 
 a totally ordered  countable group.

The restriction $L|_\Lambda$ defines
a $\Gamma$-grading on $\fl_V$ 
as follows. 
For $a\in \Gamma$, 
we set 
\[\fl_a={\rm span} \{u\in\fl_V\,|\,{\ad}_\delta(u)
=\lambda u\quad\text{for some}\quad 
\lambda\in\Lambda
\quad\text{with}\quad L(\lambda)=a\}.\]
Then 
$\fl_V=\bigoplus_{a\in \Gamma} \fl_a$.
For $\lambda,\mu\in\Lambda$, let 
$a=L(\lambda)$ and 
$b=L(\mu)$. 
We claim that $[\fl_a,\fl_b]\subset\fl_{a+b}$.
Indeed, for any $u\in\fl_a$ and any $v\in\fl_b$
such that ${\ad}_\delta(u)=\lambda u$ and 
${\ad}_\delta(v)=\mu v$, we obtain
\[{\ad}_\delta([u,v])=(\lambda+\mu)[u,v].\]
With respect to this grading, 
$\p_{\lambda_{\max}}\in V$ 
represents the principal part of $\p$, that is,
 the maximal graded component 
of $\p$. 
Indeed, 
$L(\lambda_{\max})>L(\lambda)$
for any other spectral value 
$\lambda\in\sigma$.

Let us recall that $\lambda_{\max}\neq 0$ and 
$L(\lambda_{\max})>0$.
Since 
$\p$ is locally finite of nonzero degree,
$\p_{\lambda_{\max}}\in\fh$ 
is locally nilpotent,
see  Lemma \ref{lem:vdV} 
and Remark \ref{rem:tot-ord}. 
\eproof
\bproof[Proof of Theorem \ref{thm:rk-1}]
If $\fh$ contains a nonzero locally nilpotent 
derivation, then it is triangulable;
see Corollary \ref{cor:1-lnd}. 
Otherwise, by
Lemma \ref{lem:1-lnd}, $\fh$ is 
generated by commuting 
semisimple derivations. Consequently, 
all derivations $\p\in\fh$ are semisimple;
$\fh$ is thus a toral subalgebra 
of rank at most two
and, therefore, is 
conjugate to a subalgebra of 
$\ft_2\subset\fj^+_2$
by virtue of the Gutwirth
 and Bia\l{}ynicki-Birula
linearization theorems, 
see Remark \ref{rem:ss}.
The theorem follows. 
\eproof
This completes the proof of 
Theorem \ref{main-thm}.
\subsection*{\bf Acknowledgments.}
We thank Hanspeter Kraft, Daniel Daigle,
and an anonymous referee
who pointed out errors
in earlier versions of this article.
We also thank 
 Ivan Arzhantsev  and 
 Alexander Perepechko for their interest and
helpful remarks,  Alexander Perepechko 
for his assistance in 
correcting typos and verifying proofs
by using an AI platform, 
Andriy Regeta 
for his question, and Daniel Daigle 
for stimulating discussions 
and suggestions that helped
improve the presentation. 
Our thanks also go to the anonymous 
referees for their careful reading and 
for their pertinent, constructive criticism. 
\medskip

\section{Added  in proofs.}
After this 
paper was submitted, the arxiv preprint 
\cite{Bel26} appeared. Theorem D of this 
preprint answers Question \ref{ques-1} 
in affirmative, confirming our 
conjecture. 
This result is proven in \cite{Bel26} 
for arbitrary (not necesasarily reduced) 
affine algebra 
$A$ defined over a 
field of characteristic zero 
(not necesasarily algebraically closed). 
The proofs in \cite{Bel26} are original 
and different from our ones. 

Combining our results in Section 2 and 
results in  \cite{Bel26}, we come to 
the positive answer to 
our Question \ref{ques-2}.

\bcor\la{cor:new} Let 
$h\subset \Der(\cO(X))$ 
be a solvable Lie algebra
generated by a finite set of 
locally finite derivations
and $\tilde h$ be 
as in \eqref{eq:G} 
(see Question \ref{ques-2}).
Then $\tilde h=\Lie(G)$, 
where $G\subset\Aut(X)$ is
a connected 
algebraic subgroup
as in \eqref{eq:G}. 
In particular, $\tilde h$ 
is finite-dimensional. 
\ecor
\bproof By Corollary \ref{cor:alg-gener},
$\tilde h$ is solvable and generated 
by a finite set of 
algebraic Lie subalgebras 
of $\Lie(\Aut(X))$. 
By \cite[Theorem D]{Bel26}, 
$\tilde h$ is 
finite-dimensional. 
Then \cite[Theorem A]{KZ24}
applies and gives the result 
(cf. also Remark \ref{rem:G} and 
\cite[Corollary C]{Bel26}). 
\eproof

\end{document}